\magnification =1200
\font\medtenrm=cmr10 scaled\magstep2
\font\smalletters=cmr8 at 10truept
\def\sqr#1#2{{\vcenter{\vbox{\hrule height.#2pt\hbox{\vrule width.#2pt
height#1pt\kern#1pt \vrule width.#2pt}\hrule height.#2pt}}}}
\def\square{\mathchoice\sqr64\sqr64\sqr{2.1}3\sqr{1.5}3}
\centerline {\medtenrm Linear
systems and determinantal random point fields}\par
\vskip.1in
\centerline {\medtenrm Gordon Blower}\par
\centerline {\sl Department of Mathematics and Statistics, 
Lancaster University}\par
\centerline {\sl Lancaster LA1 4YF, England UK. E-mail:
g.blower@lancaster.ac.uk}\par
\centerline {\bf 8th August 2008}\par
\vskip.05in

\hrule

\vskip.1in
\noindent {\bf Abstract}\par
\indent Tracy and Widom showed that fundamentally important kernels
in random matrix theory arise from systems of differential equations
with rational coefficients. More generally, this paper considers symmetric
Hamiltonian systems and determines the properties of kernels that
arise from them. The inverse spectral problem for 
self-adjoint Hankel operators gives sufficient conditions for a 
self-adjoint operator to be the Hankel operator on $L^2(0, \infty )$ from a 
linear system in continuous time; thus this paper expresses certain kernels as squares of Hankel 
operators. For suitable linear systems $(-A,B,C)$ with one 
dimensional input and output spaces, there exists a Hankel operator 
$\Gamma$ with kernel $\phi_{(x)}(s+t)=Ce^{-(2x+s+t)A}B$ such that 
$g_x(z)=\det (I+(z-1)\Gamma\Gamma^\dagger )$ is 
the generating function of a determinantal random point 
field on $(0,\infty )$. The inverse scattering transform 
for the Zakharov--Shabat system involves a Gelfand--Levitan 
integral equation such that the trace of the diagonal of the solution 
gives ${{\partial}\over{\partial x}}\log g_x(z)$. 
Some determinantal point fields in random matrix 
theory satisfy similar results.\par

\par
\vskip.05in
\hrule
\vskip.1in

\noindent {\bf 1. Introduction}\par
\vskip.05in

\indent Traditionally, one begins random matrix theory by defining
families of self-adjoint $n\times n$ matrices endowed with 
probability measures, known as
ensembles, and then one determines the joint distribution of the 
random
eigenvalues. By scaling the variables and letting $n\rightarrow\infty$,
 one obtains various kernels which reflect the properties of large 
random matrices. The kernels generate determinantal random point fields
 in Soshnikov's sense [16, 20]. It turns out that many such kernels in
 random matrix theory have the form 
$$K(x,y)={{f(x)g(y)-f(y)g(x)}\over {x-y}}\qquad (x,y>0)\eqno(1.1)$$
\noindent where $f$ and $g$ satisfy the system of
differential equations 
$$m(x){{d}\over{dx}}\left[\matrix{ f(x)\cr g(x)\cr}\right] 
=\left[\matrix{ \alpha (x)&\beta (x) \cr -\gamma (x)&-\alpha
(x)\cr}\right] \left[\matrix{ f(x)\cr g(x)\cr}\right] ,\eqno(1.2)$$
\noindent where $m, \alpha , \beta $ and $\gamma$ are real
polynomials. Tracy and Widom [17, 18, 19] began what amounts to a
classification of kernels that arise from such
differential equations, and their analysis revealed
detailed results about the fundamental ensembles. \par 
\vskip.05in
\noindent --------------------\par
\vskip.05in
\noindent {\smalletters This work was partially supported by EU Network
Grant MRTN-CT-2004-511953 `Phenomena in High Dimensions'.}\par
\vfill
\eject

\indent Of particular interest is the Airy kernel
$$K_\lambda (x,y)={{{\hbox{Ai}}(x-\lambda ){\hbox{Ai}}'(y-\lambda  )-{\hbox{Ai}}'(x-\lambda ){\hbox{Ai}}(y-\lambda )}\over{x-y}}\eqno(1.3)$$
\noindent on $L^2(0, \infty )$ where Airy's function ${\hbox{Ai}}$ satisfies ${\hbox{Ai}}''(x)=x{\hbox{Ai}}(x)$. Some of the fundamental properties of this kernel involve the remarkable formula
$$K_\lambda (x,y)=\int_0^\infty {\hbox{Ai}}(x+u-\lambda ){\hbox{Ai}}(u+y-\lambda )\, du\eqno(1.4)$$
\noindent which expresses the operator $K$ as the square of the Hankel operator on $L^2(0, \infty )$ that has kernel 
${\hbox{Ai}}(x+y-\lambda )$. \par
\indent The differential equation (1.2) is an example of a symmetric Hamiltonian system, as we 
can define more generally.\par
\vskip.05in
\noindent {\bf Definition} {\sl (Symmetric Hamiltonian system).} For an integer $m\geq 1$, let $J$ be the matrix
$$J=\left[\matrix{0&-I_m\cr I_m& 0\cr}\right],\eqno(1.5)$$
\noindent which satisfies $J^2=-I_{2m}$ and $J^T=-J$, and let $E(x)$ 
and $F(x)$ be $(2m)\times (2m)$ real
symmetric matrices
for each $x>0$ such that $x\mapsto E(x)$ and $x\mapsto F(x)$ are
continuous. Then we consider the symmetric Hamiltonian system
$$J{{d}\over{dx}}\Psi_\lambda =\bigl( \lambda E(x)+F(x)\bigr)
\Psi_\lambda \eqno(1.6)$$
\noindent where $\Psi_\lambda (x)$ is a $(2m)\times 1$ complex vector.
In particular, when $E(x)$ and $F(x)$ have entries that are rational functions of $x$,
we have a system considered by Tracy and Widom [19].\par
\indent Given a solution $\Psi_\lambda \in L^\infty ((0, \infty ); {\bf C}^{2m})$, we introduce the kernel
$$K_{s, \lambda } (x,y)={{\langle J\Psi_\lambda (x+s), \Psi_\lambda (y+s)\rangle_{{\bf R}^{2m}}}\over{x-y}}\qquad (x,y>0)\eqno(1.7)$$
\noindent and we investigate the properties of $K_{s,\lambda}$.\par 
\indent More generally, we consider families of kernels $K_{t, \lambda }(x,y)$ for $t, \lambda >0$, that satisfy some of the following properties as operators on $H=L^2(0, \infty )$: \par
\noindent $(1^o)$ the Lyapunov equation holds
$${{\partial}\over{\partial s}}K_{s,\lambda }=-AK_{s, \lambda
  }-K_{s, \lambda }A^\dagger\eqno(1.8)$$
\noindent as a sesquilinear form on $D(A^\dagger )\times D(A^\dagger
)$, where $e^{-sA}$ is a bounded $C_0$ semigroup on $H$ and $D(A^\dagger )$ is the domain of $A^\dagger$;\par 
\noindent $(2^o)$  $0\leq K_{s,\lambda }\leq I$ for all $s\geq s_0$ for some $s_0<\infty $;\par
\noindent $(3^o)$ $s\mapsto K_{s,\lambda }$ is decreasing and converges
strongly to $0$ as $s\rightarrow\infty$;\par
\noindent $(4^o)$ $K_{s, \lambda }$ is of trace class;\par
\noindent $(5^o)$ the operator on $H$ with kernel ${{\partial}\over{\partial s}}K_{s,\lambda }$ has finite rank.\par 
\vskip.05in

\vskip.05in
\indent In fact, many of the properties of ensembles which arise in 
random matrix theory are essentially consequences of the properties 
$(1^o)-(5^0)$, in a sense which we now make more precise. 
We recall from [16] the notion of a determinantal random point field.\par
\vskip.05in
\noindent {\bf Definition} {\sl (Configurations).} A configuration on
${\bf R}$ is an ordered list $\lambda =(\lambda_j)_{j=-\infty}^\infty$
such that $\lambda_j\leq\lambda_{j+1}$ for all $j\in {\bf Z}$; the configuration is
locally finite if 
$\nu_\lambda (L)=\sharp\{j: \lambda_j\in L\}$ is finite for all compact sets $L$. Let $\Lambda$ be the space of all locally finite configurations on ${\bf R}$. For each bounded and Borel set $E$, and $k=0, 1, \dots , $ we let
$$C^E_k=\{ \lambda\in \Lambda : \nu_\lambda (E) =k\}$$
\noindent be the set of all locally finite configurations that have
$k$ terms in $E$; now let $B$ be the $\sigma$-algebra generated by the
$C_k^E$. A random point field $({\bf P}, \Lambda , B)$ on ${\bf R}$ is a probability
measure ${\bf P}:B\rightarrow [0,1]$. We let $\nu (a,b)$ be the random
variable that gives the number of points in $(a,b)$, so $\nu (x,\infty )=\sharp\{ j: \lambda_j>x\}$.\par
\vskip.05in
\noindent {\bf Definition} {\sl (Correlation functions).} Given nonnegative integers $n_j$ such that $\sum_{j=1}^k n_j=n$ and disjoint Borel sets $E_j$ we consider $\lambda\in \Lambda$ such that
$\nu_\lambda (E_j)\geq n_j$ for all $j$. Then 
$$N_{E_j, n_j; j=1,\dots ,k} =\prod_{j=1}^k {{\nu_\lambda (E_j)!}\over 
{(\nu_\lambda (E_j)-n_j)!}}\eqno(1.9)$$
\noindent gives the number of ways of choosing 
$n_j$ points $\lambda_\ell$ from the $\nu_\lambda (E_j)$ points of 
$\lambda$ that are in $E_j$ for all $j$. The correlation function 
$R_n:{\bf R}^n\rightarrow {\bf R}_+$ for ${\bf P}$ is a locally integrable function, which
is symmetrical with respect to permutation of its variables, such that 
 $${\bf E}N_{E^{n_j}_j; j=1,\dots ,k}= \int_{E^{n_1}_1}\dots \int_{E^{n_k}_k}R_n(x_1, \dots , x_n)dx_1\dots dx_n\eqno(1.10)$$
 \noindent for all disjoint Borel sets $E_j$ $(j=1, \dots , k)$. This is the 
expected number of configurations that have $\nu_\lambda (E_j)\geq n_k$ for all $j$.\par
\vskip.05in
\indent Conversely, Soshnikov [16] observed that one
can introduce a random point
field from the determinants of kernels that satisfy minimal conditions. We state without proof the following general existence theorem for
determinantal random point fields. \par
\vskip.05in

\noindent {\bf Lemma 1.1.} {\sl Suppose that 
$K:{\bf R}\times {\bf R}\rightarrow {\bf C}$ is a
continuous kernel such that\par
\noindent $(1^o)$ the integral operator with kernel $K(x,y)$ on $L^2({\bf R})$ satisfies
$0\leq K\leq I$;\par
\noindent $(4^o)$ the kernel ${\bf I}_{[a,b]}(x)K(x,y){\bf I}_{[a,b]
}(y)$ on $L^2({\bf R})$ is of trace class for all finite $a,b$. \par
\noindent Then there exists a random point field such that the
correlation functions satisfy
$$R_n(x_1, \dots ,x_n)=\det\bigl[ K(x_j,
x_k)\bigr]_{j,k=1}^n\qquad (n=1, 2,\dots ).\eqno(1.11)$$
\indent  Further, the $R_n$ $(n=1, 2, \dots )$ uniquely determine ${\bf P}$.}\par
\vskip.05in
\noindent {\bf Definition} {\sl (Determinantal point field)}. 
If the $R_n$ have the form (1.11), then $({\bf P},\Lambda , B)$ is a determinantal random point field.\par
\vskip.05in
\indent In this paper we introduce natural examples of kernels $K$ by means of linear systems, and recover properties $(1^o)-(5^o)$ in a systematic manner. We
summarise the construction in this introduction, and describe the details in section 2. Consider an  operator $A$ with domain $D(A)$ in state space $H$ such that the $C_0$ semigroup $e^{-tA}$ is bounded, so $\Vert e^{-tA}\Vert \leq M$ for some $M<\infty$ and all $t>0$. Consider the linear system
$$\eqalignno{{{dX}\over{dt}}&=-AX+BU\qquad (X(0)=0)\cr
Y&=CX&(1.12)\cr}$$
\noindent where $B:{\bf C}\rightarrow D(A)$ is bounded and $C:D(A)\rightarrow {\bf C}$ is admissible for $A^\dagger $, so
$Y\in L^2(0,\infty )$. We let $\phi (x)=Ce^{-xA}B$ and $\phi_{(x)}(y)=\phi (y+2x)$, then introduce the Hankel operators 
$$\Gamma_{\phi_{(x)}}f(y)=\int_0^\infty \phi_{(x)}(y+u)f(u)\, 
du\eqno(1.13)$$
\noindent from a suitable domain in $L^2(0, \infty )$ into 
$L^2(0, \infty )$. We also consider the Gelfand--Levitan integral equation
$$T(x,y)-\Phi (x+y)-\int_x^\infty T(x,z)\Phi (z+y)\, dz=0\qquad (0<x\leq y<\infty )\,\eqno(1.14)$$
\noindent where $S$ and $\Phi$ are both either (i) real scalars, (ii) $2\times 2$ real diagonal matrices, or (iii) $2\times 2$ complex matrices. We state our main theorem as follows.\par

\vskip.05in
\noindent {\bf Theorem 1.2.} {\sl  Suppose that the controllability Gramian 
$$L_x=\int_x^\infty e^{-tA}BB^\dagger e^{-tA^\dagger }\, dt\qquad (x>0)\eqno(1.15)$$
\noindent is of trace class on $H$ and of operator norm $\Vert L_x\Vert <1;$ likewise suppose that the observability Gramian 
$$Q_x=\int_x^\infty e^{-tA^\dagger}C^\dagger Ce^{-tA}\,
dt\qquad (x>0)\eqno(1.16)$$
\noindent is of trace class on $H$ and that $\Vert Q_x\Vert <1.$ \par
\noindent (i) If $C=B^\dagger$ and $A=A^\dagger$, then let $g_x(z)=\det (I+(z-1)\Gamma_{\phi_{(x)}})$ and $\Phi (x)=\phi (x).$\par
\noindent (ii) If $\phi (x)$ is real, then let $g_x(z)=
\det (I+(z-1)\Gamma_{\phi_{(x)}}^2)$ and} $\Phi (x)={\hbox{diag}}[-\phi (x), \phi (x)].$ \par
\noindent {\sl (iii) Or let $g_x(z)=\det (I+(z-1)\Gamma_{\phi_{(x)}}
\Gamma_{\phi_{(x)}}^\dagger)$ and 
$$\Phi (x)=\left[\matrix{0& \bar \phi (x)\cr -\phi
(x)&0\cr}\right]\qquad (x>0).$$
\noindent Then in each case there exists a determinantal random point field on $(0, \infty )$ with generating function $g_x(z)={\bf E}z^{\nu (x, \infty )}$ such that}
$${{\partial}\over{\partial x}}\log g_x(0)=
{\hbox{trace}}\, T(x,x)\qquad (x>0)\eqno(1.17)$$
\noindent {\sl is given by the diagonal of the solution of the 
Gelfand--Levitan integral equation (1.14).}\par
\vskip.05in  
\indent The integral equation in case (i) is associated with the inverse scattering problem for the Schr\"odinger equation on the real line and in case (ii) by a pair of Schr\"odinger equations; whereas the integral equation in (iii) is associated with a Zakharov--Shabat system.\par

\noindent The fundamental examples of determinantal random point fields in random matrix theory involve kernels associated with self-adjoint Hamiltonian systems of
differential equations. 
In section 3 we introduce the notion of a symmetric Hamiltonian system
with matrix potential, as considered previously by Atkinson and many
others; see [5]. We consider spatial kernels 
$K_\lambda$ associated with symmetric Hamiltonian systems, and 
give a sufficient condition for the kernel to factor as 
$K_\lambda =\Gamma_\lambda^\dagger \Gamma_\lambda$, where 
$\Gamma_\lambda$ is a vectorial Hankel operator. 
As we show in section 4, this covers some fundamental examples 
of kernels that arise in random matrix theory, and we recover case (ii) of Theorem 1.2. A similar computation shows how (iii) arises.\par 
\indent Schr\"odinger differential operators on 
$L^2(0, \infty )$ with bounded potentials give rise to kernels 
in the spectral variables which satisfy
$(5^o)$, as we discuss in section
5. The Korteweg--de Vries flow has a natural effect on the kernels.
 In section 6, we consider the Zakharov--Shabat system and 
establish case (iii) of Theorem 1.2; here the kernels behave naturally
under the flow associated with the  cubic nonlinear Schr\"odinger
equation. 
Some of the calculations will be familiar to specialists in 
the theory of scattering from [1, 6, 23], but we include them here so that the 
paper is self-contained.\par

\vskip.05in

\noindent {\bf 2. Linear systems and their determinants}\par

\vskip.05in
\noindent {\bf Definition} {\sl (Linear system).} Let $H$ be a separable complex Hilbert space, called the state space, and $H_0$ a separable complex Hilbert space called the output space. Let $e^{-sA}$ be a $C_0$ semigroup on $H$, such that $\Vert e^{-sA}\Vert\leq M$ for some $M<\infty$ and all $s>0$, and let $D(A)$ be the domain of the generator $-A$, which is a dense linear subspace of $H$ and itself a Hilbert space for the norm $\Vert\xi\Vert_{D(A)}=\bigl( \Vert\xi\Vert^2_H+\Vert A\xi \Vert^2_H\bigr)^{1/2}$.
In the language of linear systems from [14, 15], we consider the continuous time system
$$\eqalignno{{{dX}\over{dt}}&=-AX+BU\qquad (t>0)\cr
Y&=CX,\qquad X(0)=0.&(2.1)}$$
\noindent where $B:H_0\rightarrow D(A)$ and $C:D(A)\rightarrow H_0$ are 
bounded linear operators; this is known as a $(-A, B, C)$ system. 
Let $\phi (x)=Ce^{-Ax}B$, so $\phi\in L^\infty ((0, \infty ); 
{\hbox{B}}(H_0)).$ The associated Hankel operator $\Gamma_\phi$ is the integral operator 
$$\Gamma_\phi f(x)=\int_0^\infty \phi (x+y)f(y)\, dy,\eqno(2.2)$$
\noindent defined from some dense linear subspace of 
$L^2((0, \infty ); H_0)$ into $L^2((0,\infty );H_0)$.\par
\indent We introduce the transfer function
$$\hat\phi (\lambda )=C(\lambda I+A)^{-1}B,\eqno(2.3)$$
\noindent which we recognise as the Laplace transform of 
$\phi (x)=Ce^{-Ax}B$; the Fourier transform of $\phi$ gives the
scattering data.
\noindent Suppose that $U\in L^2((0, \infty); H)$ has Laplace
transform $U(\lambda )$, and that $\hat\phi :{\bf C}_+\rightarrow {\hbox{B}}(H_0)$ is a bounded analytic function. Then $Y\in L^2((0, \infty
);H_0)$ has  Laplace transform $\hat Y$ such that 
$\hat Y(\lambda )=\hat\phi (\lambda )\hat U(\lambda ).$\par
\vskip.05in

\noindent {\bf Definition} {\sl (Admissible).}  We say 
that a bounded linear operator $C:D(A)\rightarrow H_0$ is admissible for $e^{-sA}$ if $Ce^{-sA}\xi$ belongs to $L^2((0, \infty ); H_0)$ for all $\xi\in H$, and there exists $K_C(A)$ such that    
$$\int_0^\infty \Vert Ce^{-sA}\xi\Vert^2_{H_0}\, ds\leq K_C(A)^2\Vert \xi\Vert_H^2\qquad (\xi\in H),\eqno(2.4)$$
\noindent equivalently, the operator $\Theta^\dagger :H\rightarrow
L^2((0, \infty ); H_0)$ is bounded where 
$\Theta^\dagger \xi =Ce^{-sA}\xi$ and $\Vert \Theta\Vert =K_C(A)$. 
Examples in [9] show that the notion of admissibility is difficult to characterize simply.  

\vskip.05in 

\noindent {\bf Definition} {\sl (Schatten ideals).} Let $c^2$ be the space of Hilbert--Schmidt operators, and $c^1$ be the space of trace class operators, with the usual norms.\par


\vskip.05in
\noindent {\bf Proposition 2.1.} {\sl Suppose that $C$ is admissible
  for $e^{-sA}$ and that $B:H_0\rightarrow D(A^\dagger )$ has
  $B^\dagger$ admissible for $e^{-sA^\dagger}$. Then the observability Gramian 
$$Q_x=\int_x^\infty e^{-sA^\dagger}C^\dagger Ce^{-sA}\, ds\qquad (x>0)\eqno(2.5)$$
\noindent and the controllability Gramian
$$L_x=\int_x^\infty e^{-sA}BB^\dagger e^{-sA^\dagger}\, ds\qquad (x>0)\eqno(2.6)$$
\noindent define bounded linear operators $Q_x, L_x:H\rightarrow H$ by
these strongly convergent integrals such that:\par
\noindent $(1^o)$ the derivatives satisfy the Lyapunov equations 
$${{\partial Q_x}\over{\partial x}}=-A^\dagger Q_x-Q_xA, \quad 
{{\partial L_x}\over{\partial x}}=-AL_x-L_xA^\dagger\eqno(2.7)$$  
\noindent as sesquilinear forms on $D(A)\times D(A)$ and $D(A^\dagger)\times D(A^\dagger)$ respectively;\par
\noindent $(2^0)$ $0\leq Q_x  \leq K_C(A)^2I$ and $0\leq L_x\leq K_{B^\dagger}(A^\dagger)^2I$ for all $x\geq 0$;\par
\noindent $(3^0)$ $Q_x$ and $L_x$ decrease strongly to zero as $x$ increases to infinity.}\par
\noindent  {\sl $(4^o)$ Suppose further that $C(iyI+A)^{-1}$ is a Hilbert--Schmidt operator for all $y\in {\bf R}$ and that $\int_{-\infty}^\infty \Vert C(iyI+A)^{-1}\Vert^2_{c^2}dy<\infty$. Then $Q_x$ is trace class for each $x>0$, and} 
$${\hbox{trace}}\, Q_0={{1}\over{2\pi}}\int_{-\infty}^\infty \Vert C(ixI+A)^{-1}\Vert_{c^2}^2\, dx.\eqno(2.8)$$
\noindent {\sl $(5^o)$ Suppose that $H_0$ has finite dimension $m$.
Then} ${\hbox{rank}}{{\partial Q_x}\over{\partial x}}\leq m$.\par  
\vskip.05in
\noindent {\bf Proof.} $(2^0)$,$(3^0)$ The integrals converge by the definition of admissibility, and the other statements are immediate consequences.\par
\indent $(1^o)$ For $\xi\in D(A^\dagger)$, the $e^{-sA^\dagger}\xi$ is differentiable in $H$  with derivative $-e^{-sA^\dagger}A^\dagger\xi$. By the fundamental theorem of calculus, we have 
$$-AL_x-L_xA^\dagger=\int_x^\infty{{d}\over{ds}}\bigl(e^{-sA }
BB^\dagger e^{-sA^\dagger}\bigr)\, ds=-e^{-xA}BB^\dagger e^{-xA^\dagger}
,\eqno(2.9)$$
\noindent as bilinear forms on $D(A^\dagger)\times D(A^\dagger)$,
 hence the result.\par 
\indent  $(4^0)$ Let $(e_j)_{j=1}^{\infty}$ be an orthonormal basis for $H$. By Plancherel's formula in Hilbert space, we have
$$\int_0^\infty\Vert Ce^{-yA}e_j\Vert_{H_0}^2dy={{1}\over{2\pi}}\int_{-\infty}^\infty \Vert C(ixI+A)^{-1}e_j\Vert_{H_0}^2\, dx\eqno(2.10)$$
\noindent and summing this identity we deduce
$$\sum_{j=1}^\infty \langle Q_0e_j, e_j\rangle ={{1}\over{2\pi}}\int_{-\infty}^\infty \sum_{j=1}^\infty \Vert C(ixI+A)^{-1}e_j\Vert_{H_0}^2\, dx\eqno(2.11)$$
\noindent and hence
$${\hbox{trace}}\, Q_0=\Vert Q_0^{1/2}\Vert^2_{c^2}={{1}\over{2\pi}}\int_{-\infty}^\infty \Vert C(ixI+A)^{-1}\Vert^2_{c^2}\, dx,\eqno(2.12)$$
\noindent so $Q_0$ is trace class.\par
\indent $(5^0)$ From the expression
$${{\partial Q_x}\over{\partial x}}=-e^{-xA^\dagger}C^\dagger
Ce^{-xA}\eqno(2.13)$$
\noindent it follows that the rank of ${{\partial Q_x}\over{\partial
    x}}$ is less than or equal to the rank of $C$, hence is less than or
    equal to $m$.\par  
\rightline{$\square$}\par
\vskip.05in

\noindent {\bf Proposition 2.2} 
(Determinant of the observability Gramian).\par
\noindent {\sl Suppose that $(4^0)$ holds, so that the observability operator 
$\Theta :L^2((0, \infty );H_0)\rightarrow H$ is Hilbert--Schmidt, where
$$\Theta f=\int_0^\infty e^{-sA^\dagger} C^\dagger f(s)\,ds
\qquad (f\in L^2((0, \infty ); H_0)).\eqno(2.14)$$
\noindent (i) Then  
$$\det (I-\lambda Q_x)=\det (I-\lambda P_{(x, \infty )} \Theta^\dagger \Theta P_{(x, \infty )})\qquad (\lambda \in {\bf C}, x\geq 0).\eqno(2.15)$$
\noindent defines an entire function that has all its zeros on the positive real axis.\par 
\noindent (ii) Suppose further that $A=A^\dagger$. Then the Hankel operator $\Gamma_\phi$ on 
$L^2((0, \infty ); H_0)$ with kernel $\phi (s+t)=Ce^{-(s+t)A}C^\dagger$, has $\Gamma_\phi=\Theta^\dagger\Theta\geq 0$ and}
$${{\partial }\over{\partial x}}{\hbox{trace}} \,
Q_x=-{\hbox{trace}}\, \phi (2x ).$$
\noindent {\sl (iii) Suppose still further that $H_0={\bf C}^m$ where $m<\infty$. 
Then the zeros of $\det (I-\lambda Q_x)$ have order less than or equal to $m$.}\par
\vskip.05in
\noindent {\bf Proof.} (i) We have $\Theta^\dagger \xi (t)=Ce^{-tA}\xi $ and hence $\Theta\Theta^\dagger =Q_0$. Further, since the operators $\Theta$ and $\Theta^\dagger$ are Hilbert--Schmidt, we can rearrange terms in the determinant and obtain
$$\det (I-\lambda P_{(x, \infty )} \Theta^\dagger \Theta P_{(x, \infty )})=\det (I-\lambda \Theta P_{(x, \infty )} \Theta^\dagger )=\det (I-\lambda Q_x).\eqno(2.16)$$
\noindent The zeros of $\det (I-\lambda Q_x)$ are $1/\lambda_j$, where $\lambda_j$ are the positive eigenvalues of $Q_x$.\par 
\indent (ii)  Now
$$ \Theta^\dagger \Theta f(t)=Ce^{-tA}\int_0^\infty e^{-sA^\dagger} 
C^\dagger f(s)\, ds,\eqno(2.17)$$
\noindent so $\Theta^\dagger\Theta$ reduces to a Hankel operator 
when $A=A^\dagger .$  Further, we have
$$\eqalignno{{\hbox{trace}}\,Q_x&=
\int_x^\infty {\hbox{trace}}\, e^{-tA}C^\dagger Ce^{-tA}\, dt\cr
&=\int_x^\infty {\hbox{trace}}\, Ce^{-2tA}C^\dagger \, dt,&(2.18)}$$
\noindent whence the result.\par
\indent (iii) The (block) Hankel operator with kernel
 $\phi_{(x)}(s+t)=Ce^{(s+t+2x)A}C^\dagger$ is non negative and compact, 
and is unitarily equivalent to the some matrix 
$[a_{j+k}]_{j,k=1}^\infty$ which is made up of $m\times m$ blocks. 
Hence its spectrum consists of $0$ together with a sequence of 
eigenvalues $\lambda_j$ of multiplicity less than or equal to $m$
 which decrease strictly to $0$ as $j\rightarrow\infty$ by 
[14, Theorem 2]. Hence the zeros of the function $\det (I-\lambda P_{(x, \infty )} \Theta^\dagger \Theta P_{(x, \infty )})$ have order less than or equal to $m$ at the points $1/\lambda_j$.\par
\rightline{$\square$}\par
\vskip.05in
\indent To express the hypotheses of Proposition 2.2(ii) in terms of spectra, we present the following result, which is known to specialists.\par
\vskip.05in
\noindent {\bf Definition} {\sl (Carleson measure).} Let $\mu$ be a positive Radon measure on 
${\bf C}_+=\{ z\in {\bf C}: \Re z>0\}$. Then $\mu$ is a Carleson measure if there exists $c_0>0$ such that 
$$\mu \bigl( [0,x]\times [y-x,y+x]\bigr)\leq c_0x\qquad (x>0, y\in {\bf R}).\eqno(2.19)$$
\vskip.05in
\noindent {\bf Proposition 2.3.} {\sl Suppose that $A$ is 
self-adjoint and has purely discrete spectrum $(\kappa_j)$, with $\kappa_j>0$ listed according to multiplicity, and that $(e_j)$ is  corresponding orthonormal basis of eigenvectors. Let $\phi (x)=Ce^{-xA}C^\dagger$.\par
\noindent (i) Then $\Gamma_\phi$ is bounded if and only if $\sum_{j=1}^\infty \vert Ce_j\vert^2\delta_{\kappa_j}$ is a Carleson measure.\par
\noindent (ii) If $\sum_{j=1}^\infty \vert Ce_j\vert^2/\kappa_j$ converges, then $\Gamma_\phi$ is of trace class.}\par
\vskip.05in
\noindent {\bf Proof.} (i) We use hats to denote Laplace transforms, 
and let $H^2$ be the usual Hardy space on ${\bf C}_+$ as in [10]. By the Paley--Wiener theorem, the Laplace transform gives a unitary map from $L^2(0, \infty )$ onto $H^2$. Then 
$$\eqalignno{\langle \Gamma_\phi f,f\rangle
&=\int_0^\infty\int_0^\infty \sum_{j=1}^\infty \vert 
Ce_j\vert^2 e^{-(s+t)\kappa_j}f(s)\overline{f(t)}dsdt\cr
&= \sum_{j=1}^\infty \vert Ce_j\vert^2 \vert \hat
f(\kappa_j)\vert^2.&(2.20)}$$
\noindent Hence $\Gamma_\phi$ is bounded if and only if  there exists $c_1$ such that 
$$\Bigl\langle \sum_{j=1}^\infty \vert Ce_j\vert^2
\delta_{\kappa_j},  \vert \hat f\vert^2\Bigr\rangle \leq 
c_1\lim_{x\rightarrow 0+}\int_{-\infty}^\infty \vert \hat f(x+iy)
\vert^2dy;\eqno(2.21)$$
\noindent which holds if and only if we have a Carleson measure;
see [10].\par
\indent (ii) Note that $\sqrt{2\kappa_j}/(z+\kappa_j)$ is a unit vector
 in $H^2,$ so $\hat f\mapsto 2\kappa_j \hat f(\kappa_j)/(
z+\kappa_j)$ has rank one and norm one as an operator on $H^2$; hence the result by convexity.\par  
\rightline{$\square$}\par

\noindent {\bf Definition} {\sl (Balanced system).} If $Q_0=K_0$, then the system is balanced.\par
\vskip.05in
\noindent {\bf Remarks}\par 
\noindent (i) The controllability operator $\Xi :L^2((0, \infty ); H_0)\rightarrow H$
$$\Xi f=\int_0^\infty e^{-tA}B f(t)\, dt\eqno(2.22)$$
\noindent satisfies an obvious analogue of Proposition 2.2. Note that $\Gamma_\phi=\Theta^\dagger\Xi .$\par
\noindent (ii) One can interchange the controllability and 
observability operators by interchanging 
$(-A, B, C)\leftrightarrow (-A^\dagger ,C^\dagger, B^\dagger )$, which interchanges $\Gamma_\phi\leftrightarrow\Gamma_\phi^\dagger.$\par
\noindent However, we will consider in section 5 some self-adjoint Hankel operators which do not arise from the special case of $A=A^\dagger$ and $B=C^\dagger$.\par
\noindent (iii) The set ${\bf K}$ of kernels that satisfy $(2^o), (3^0), (4^o)$ and $(5^o)$
is convex; further, for $K\in {\bf K}$ and $U\in B(H)$ such that 
$\Vert U\Vert\leq 1$, we have $U^\dagger KU\in {\bf K}$.\par 
\noindent (iv) If $(1^o)$ holds with a finite-dimensional $H$, then
$(5^o)$ holds. But $(5^o)$ is not implied by finite-dimensionality of $H_0$.\par
\noindent (v) For $x>0$, the shifted system $(-A, e^{-xA}B, Ce^{-xA})$
has observability operator $Q_x$,
controllability operator $L_x$ and, with 
$\phi_{(x)}(t)=Ce^{-(2x+t)A}B$, the corresponding Hankel operator
is $\Gamma_{\phi_{(x)}}$.\par
\vskip.05in
\noindent {\bf Proposition 2.4.} 
(Determinants involving the Hankel operator) {\sl Suppose that the 
controllability operator $\Theta_x$ and the observability 
operator $\Xi_x$ for $(-A, e^{-xA}B, e^{-xA}C)$ are Hilbert--Schmidt. Then the operator
$R_x:H\rightarrow H$, defined by
$$R_x\xi =\int_x^\infty e^{-yA}BCe^{-yA}\xi dy,\eqno(2.23)$$ 
\noindent is of trace class and satisfies}
$$\det (I-\lambda \Gamma_{\phi_{(x)}} )=\det (I-\lambda R_x).\eqno(2.24)$$
\vskip.05in
\noindent {\bf Proof.} By Proposition 2.2, the operator $R_x$ is 
trace class. By rearranging, we obtain
$$\eqalignno{\det (I-\lambda \Gamma_{\phi_{(x)}} )&
=\det (I-\lambda \Theta_x^\dagger\Xi_x)\cr
&=\det (I-\lambda \Xi_x\Theta_x^\dagger)\cr
&=\det (I-\lambda R_x).&(2.25)}$$
\rightline{$\square$}\par
\vskip.05in
 
\vskip.05in
\indent In section 3 we introduce some kernels that arise from 
Hankel operators as in Proposition 2.2. The kernels are defined with symmetric Hamiltonian systems, as we recall in section 3.\par
\vskip.05in

\noindent {\bf 3. Kernels arising from Hamiltonian systems of ordinary differential
equations}\par 
\vskip.05in
\indent Let $D$ be a domain that is symmetrical with respect to the real axis and contains $(0, \infty )$. We later define kernels $K_\lambda$ that satisfy some of the following properties:\par
\noindent $(6^o)$ $K_\lambda$ defines a bounded linear operator on $H$ for all $\lambda\in D$;\par
\noindent $(7^o)$ $\lambda\mapsto K_\lambda$ is analytic on $D$;\par
\noindent $(8^o)$ $K_{\bar\lambda}=K_\lambda^\dagger$ for all $\lambda\in D$;\par
\noindent $(9^o)$ $K_\lambda$ is a Hilbert--Schmidt operator for all $\lambda\in D$;\par 
\noindent $(10^o)$ $K_\lambda$ is an integrable operator on $L^2(I;
dx)$, for some interval $I$; so there
exist locally bounded and measurable functions $\psi_k$ and $\xi_k$
such that 
$$K_\lambda (x,y)=\sum_{j=1}^m {{\psi_j(x;\lambda )\xi_j(y;\lambda )}
\over{x-y}}\qquad (x,y\in I; x\neq y)$$
\noindent and $\sum_{j=1}^m \psi_j(x;\lambda )\xi_j(x;\lambda )=0$.\par

\vskip.05in

\noindent {\bf Definition} {\sl (Hamiltonian system).} Let $E(x)$ 
and $F(x)$ be $(2m)\times (2m)$ real
symmetric matrices
for each $x>0$ such that $x\mapsto E(x)$ and $x\mapsto F(x)$ are
continuous. Then we consider the symmetric Hamiltonian system
$$J{{d}\over{dx}}\Psi_\lambda (x)=\bigl( \lambda E(x)+F(x)\bigr)
\Psi_\lambda (x)\eqno(3.1)$$
\noindent where $\Psi_\lambda (x)$ is a $(2m)\times 1$ complex vector. 
Suppose that for some $\lambda \in {\bf C}$, the solution
$\Psi_\lambda$ of (1.4) belongs to $L^\infty ((0, \infty ); 
{\bf C}^m).$ With the bilinear form 
$\langle (z_j), (w_j)\rangle =\sum_{j=1}^{2m}z_jw_j$, let 
$$K_\lambda (x,y) ={{\Psi_\lambda^T (y)J\Psi_\lambda (x)}\over {x-y}}=
{{\langle J\Psi_\lambda (x), \Psi_\lambda (y)\rangle}\over
{x-y}};\eqno(3.2)$$
as in l'H\^opital's rule, the diagonal of the kernel is taken to be
$$K_\lambda (x,x)=\langle \Psi_\lambda (x), (\lambda E(x)+F(x))
\Psi_\lambda (x)\rangle .$$
\vskip.05in
\noindent {\bf Proposition 3.1.} {\sl (i) Suppose that the solution
$\Psi_\lambda$ belongs to $L^\infty ((0, \infty ); {\bf C}^{2m})$ for
all $\lambda\in D$. 
Then $K_\lambda$ defines the kernel of a bounded linear operator $K_\lambda$ on $L^2((0, \infty ); {\bf
C})$ such that $K_\lambda^{\dagger}=K_{\bar\lambda },$ so $(6^o),
(7^o), (8^o)$ and 
$(10^o)$ hold.\par
\noindent (ii) Suppose that $E$ and $F$ are bounded and that 
$\Psi_\lambda $ is a solution in $L^2((0, \infty ); {\bf C}^{2m})$ of 
(3.2). Then $K_\lambda (x,y)$ defines a Hilbert--Schmidt kernel on 
$L^2((0, \infty );dx)$, so $(9^o)$ also holds.}\par
\vskip.05in
\noindent {\bf Proof.} (i) Indeed, the Hilbert transform $H$ with kernel 
$1/(\pi (x-y))$ is bounded on $L^2({\bf R})$ and $K_\lambda$ is a 
composition of $H$ and bounded multiplication operators. 
The conditions $(6^o), (7^o)$ and $(8^o)$ follow from basic facts 
about differential equations as in [8].\par 
\indent Observe that $\langle J\Psi_\lambda (x), \Psi_\lambda (x)
\rangle =0$, since we have the bilinear product; 
so the formula for $K_\lambda$ extends by continuity to a continuous
function on $(0, \infty )^2$ and $K_\lambda$ is an integrable kernel
as in $(10^o)$.\par
\indent (ii) There exist constants $c_1$ and $c_2$ such that 
$\Vert E(x)\Vert\leq c_1$ and $\Vert F(x)\Vert\leq c_2$; hence 
the differential equation gives $\Vert\Psi'_\lambda (x)\Vert\leq 
(c_1\vert \lambda \vert+c_2)\Vert\Psi_\lambda (x)\Vert .$ We deduce 
that $\Psi_\lambda' $ belongs to $L^2((0, \infty ); {\bf C}^{2m})$, 
and it is then easy to see that $\Psi_\lambda $ is bounded. A further application of the differential equation shows that $\Psi'_\lambda $ is also bounded.\par
\indent We split the Hilbert--Schmidt integral as
$$\int_0^\infty \!\!\!\int_0^\infty \vert K_\lambda (x,y)\vert^2 dxdy\leq \int\!\!\!\int_{\vert x-y\vert \geq 1}{{\Vert \Psi_\lambda (x)\Vert^2\Vert\Psi_\lambda (y)\Vert^2}\over{\vert x-y\vert^2}} dxdy$$
$$+\int\!\!\!\int_{\vert x-y\vert\leq 1}\Bigl\vert
 \Bigl\langle J\Psi_\lambda (x), {{\Psi_\lambda (x)-
\Psi_\lambda (y)}\over{x-y}}\Bigr\rangle\Bigr\vert^2 dxdy.\eqno(3.3)$$
\noindent The preceding estimates show that both of these integrals converge.\par  
\rightline{$\square$}\par
\vskip.05in  
\noindent {\bf Definition} {\sl (Shift).} For $t>0$,
 let $S_t:L^2(0, \infty )\rightarrow L^2(0, \infty )$ be the shift
 $S_tf(x)=f(x-t)$, so that $S^\dagger_tS_t=I$, and 
$S_tS^\dagger_t=P_{(t, \infty )}$ is the orthogonal projection 
onto $L^2[t, \infty )\subset L^2[0, \infty )$. In the remainder of 
this section, we are concerned with the effect of the shift on 
solutions $S_t:\Psi_\lambda (x) \mapsto \Psi_\lambda (x-t)$, and 
consequently on the kernels $K_{t,\lambda }=S_t^\dagger K_\lambda S_t.$
\vskip.05in
\noindent {\bf Definition.} For $I$ an interval in ${\bf R}$, let $L:I\times I\rightarrow M_{m}({\bf C})$ be a continuous function. We write $L\succeq 0$ if there exist continuous functions $E_j:I\rightarrow M_m({\bf C})$ for $j=1, 2,\dots $ such that $\sup_{x\in I}\Vert \sum_{j=1}^\infty E_j(x)^\dagger E_j(x)\Vert <\infty$ and 
$$L(x,y)=\sum_{j=1}^\infty E_j(y)^\dagger E_j(x)\qquad (x,y\in I).\eqno(3.4)$$
\vskip.05in
\noindent {\bf Lemma 3.2.} {\sl Let $\nu$ be probability measure 
$\nu$ on $I$, and suppose that $L\succeq 0$ on $I\times I$.\par
\noindent (i) Then $\Phi :L^2(\nu ;{\bf C}^m)
\times L^2(\nu ; {\bf C}^m)\rightarrow {\bf C}$, 
is a positive sesquilinear form, where}
$$\Phi (\xi , \eta )=\int\!\!\!\int_{I\times I} \langle L(x,y)\xi (x), \eta (y)\rangle_{{\bf C}^m} \nu (dx)\nu (dy)\qquad (\xi , \eta\in L^2(I; \nu ; {\bf C}^m)).$$
\noindent {\sl (ii) Suppose further that the defining sum (3.4) for $L$ is finite. Then $\Phi :L^2(\nu ;{\bf C})\rightarrow L^2(\nu ;{\bf C}^m)$ has finite rank.}\par 
\vskip.05in
\noindent {\bf Proof.} (i) The kernel is positive since
$$\sum_{k,\ell} a_{k,\ell}\Phi (\xi_k , \xi_\ell) =\sum_{j=1}^\infty \sum_{k,\ell}a_{k,\ell}
\Bigl\langle \int_I E_j(x)\xi_k (x)\nu (dx) , \int_I E_j(y)\xi_\ell (y)\,\nu (dy )\Bigr\rangle \geq 0.$$
\indent (ii) This is clear, since $\Phi$ can be expressed a finite tensor.\par
\rightline{$\square$}\par
\vskip.05in
\noindent {\bf Proposition 3.3.} {\sl Suppose that for some
$\lambda >0$ there exists a bounded and continuous solution $\Psi_\lambda \in L^2((0, \infty
); {\bf R}^{2m})$ to (3.2), where the coefficients $E(x)$ and $F(x)$ are bounded and satisfy  
$${{E(x)-E(y)}\over{x-y}}\succeq 0,\quad {{F(x)-F(y)}\over{x-y}}\succeq 
0\qquad (x,y\in (0, \infty )).\eqno(3.5)$$
\noindent (i) Then $K_{t, \lambda }=S_t^\dagger K_\lambda S_t$ satisfies
$(1^o), (2^o), (3^o), (6^o)$,$(7^o), (8^0), (9^o)$ and $(10^o)$.\par
\noindent (ii) The kernel $K_{t,\lambda}$ is of trace class, as in $(4^o)$  and satisfies}
$${\hbox{trace}}K_{t, \lambda }=\int_t^\infty \langle \Psi_\lambda (x), (\lambda E(x)+F(x))\Psi_\lambda (x)\rangle \, dx.$$
\noindent {\sl (iii)  If the sums involved in (3.4) for $(E(x)-E(y))/(x-y)$ and $(F(x)-F(y))/(x-y))$ are finite, then $(5^0)$ also holds.}\par
\vskip.05in
\noindent  {\bf Proof.} (i) We observe that if $L\succeq 0$, then 
$\langle L(x,y)\xi , \xi\rangle$ gives the kernel of a positive 
definite operator. $K_\lambda (z+t,w+t)$ gives the kernel that represents $S_t^\dagger K_\lambda S_t$, and hence satisfies the Lyapunov equation
$${{\partial }\over{\partial t}}K_{t, \lambda}=-AK_{t, \lambda}-K_{t, \lambda}A^\dagger$$
\noindent where $-A={{\partial}\over{\partial x}}$ generates the semigroup $S_t^\dagger$.\par
\indent  From the differential equation, we have 
$$\eqalignno{{{\partial }\over{\partial t}}K_\lambda
(x+t,y+t)=&-\lambda 
\Bigl\langle {{E(x+t)-E(y+t)}\over{x-y}}
\Psi_\lambda (x+t), \Psi_\lambda (y+t)\Bigr\rangle \cr 
&\quad -\Bigl\langle {{F(x+t)-F(y+t)}\over{x-y}}
\Psi_\lambda (x+t), \Psi_\lambda (y+t)\Bigr\rangle, &(3.6)}$$
\noindent and by the hypotheses on $E$ and $F$ we deduce that there 
exist $E_j, F_j:(0, \infty )\rightarrow M_{2m}({\bf C})$ such that 
$$\eqalignno{{{\partial }\over{\partial t}}K_\lambda
(x+t,y+t)=&-\lambda\sum_{j=1}^\infty \langle E_j(x+t)
\Psi_\lambda (x+t), E_j(y+t)\Psi_\lambda (y+t)\rangle\cr
&\quad -\sum_{j=1}^\infty \langle F_j(x+t)\Psi_\lambda (x+t), F_j(y+t)
\Psi_\lambda (y+t)\rangle .&(3.7)}$$
\noindent The right-hand side gives the kernel of a negative
 definite operator, so
$${{\partial }\over{\partial t}}\langle S_t^\dagger K_\lambda 
S_tf,f\rangle\leq 0\eqno(3.8)$$
\noindent for all $f\in L^2(0, \infty )$.\par
\indent By arguing as in Proposition 3.2, we see that $K_\lambda:L^2(0, \infty )\rightarrow L^2(0, \infty )$ is compact. For $f\in L^2(0, \infty )$, we observe that $S_tf\rightarrow 0$ weakly as $t\rightarrow\infty$ and 
since $K_\lambda $ is compact $K_\lambda S_tf\rightarrow 0$ in norm as $t\rightarrow\infty$; hence $\langle S_t^\dagger K_\lambda S_t f, f\rangle \rightarrow 0$. Consequently $\langle K_{t,\lambda} f,f\rangle$ decreases to $0$ as $t\rightarrow\infty$. Further, $S_t^\dagger K_\lambda S_t\rightarrow 0$ in Hilbert--Schmidt norm as $t\rightarrow\infty$, so there exists $s_0$ such that $\Vert S_t^\dagger K_\lambda S_t\Vert\leq 1$ for all $s\geq s_0$.\par
\indent (ii) We have proved that $K_{t, \lambda }\geq 0$, so the kernel is positive definite and continuous. By Mercer's trace formula,
$${\hbox{trace}}K_{t, \lambda}= \int_0^\infty K_{t, \lambda}(x,x)\, dx
=\int_t^\infty K_{0, \lambda}(x,x)dx.\eqno(3.9)$$
\indent (iii) If the sum over $j$ has finitely many terms, then
the expression for 
${{\partial}\over{\partial t}}K_\lambda (x+t, y+t)$ is a finite 
tensor and hence a finite-rank operator.\par  
\rightline{$\square$}\par
\vskip.05in
\noindent We now relate the notion of positivity from the previous definition to matrix monotonicity in Loewner's sense.\par

\vskip.05in
\noindent {\bf Definition} {\sl (Matrix monotone).}  Let $I$ be an open real interval and let 
$I^c={\bf R}\setminus I$. Suppose that $E :{\bf C}\setminus
I^c\rightarrow M_m({\bf C})$ is an analytic function such that
$E(x)=E(x)^{\dagger}$ for all $x\in I$ and 
$$(E(z)-E(z)^\dagger )/(2i)=\Im E(z)\geq 0\qquad (\Im z>0).\eqno(3.10)$$
\noindent Then $E$ is a Loewner matrix function on $I$; equivalently,
$E$ is said to be matrix monotone.\par
\vskip.05in
\noindent {\bf Theorem 3.4.} {\sl Suppose that for some $\varepsilon >0$
the functions $z\mapsto E(z)$ and $z\mapsto F(z)$ are matrix Loewner
functions on $(-\varepsilon , \infty )$. Suppose further that for
$\Re \lambda >0$ there exists a bounded and continuous solution $\Psi_\lambda \in L^2((0, \infty
); {\bf R}^{2m})$ to (3.1) such that $\Psi_\lambda (x)\rightarrow 0$ as
$x\rightarrow\infty$. \par
\noindent (i) Then for all $\Re \lambda >0$ there exists $\phi :(0, \infty )\rightarrow H_0$ and a bounded 
Hankel operator $\Gamma_{\phi_\lambda} :L^2(0, \infty )\rightarrow L^2((0, \infty ); H_0)$ such that 
$$K_\lambda =\Gamma_{\phi_\lambda}^\dagger \Gamma_{\phi_\lambda },
\qquad (\Re \lambda >0)\eqno(3.11)$$
\noindent and the family of kernels $K_{t, \lambda }
=S_t^\dagger K_\lambda S_t$, for $\Re \lambda >0$, satisfies conditions $(1^o)$-$(4^o)$, and $(6^o)-(8^o), (10^o)$.}\par 
\noindent {\sl (ii) If $H_0$ has finite dimension, then 
$K_{t, \lambda }$ also satisfies $(5^o)$.\par
\noindent (iii) If $H_0={\bf C}$, then $\Gamma_{\phi_\lambda}$ is a 
scalar Hankel operator and 
$K_\lambda =\Gamma_{\phi_\lambda}^\dagger \Gamma_{\phi_\lambda}$.\par
\noindent (iv) If $H_0={\bf C}$ and $\phi_\lambda$ is real-valued, 
then $K_\lambda =\Gamma^2_{\phi_\lambda}$.}\par 
\vskip.05in
\noindent {\bf Proof.} (i) We need to obtain a suitable 
$\phi_\lambda$ for the vectorial Hankel operator. Let $D ={\bf C}\setminus (-\infty, \varepsilon ]$, and let $K_\lambda (z,w)$ be a kernel on $D\times D$.
Then $K_\lambda (z+t,w+t)$ gives the kernel that represents $S_t^\dagger K_\lambda S_t$, and $K_\lambda =\Gamma_{\phi_\lambda}^\dagger \Gamma_{\phi_\lambda}$ 
if 
$$\Bigl({{\partial }\over{\partial t}}\Bigr)_{t=0}K_\lambda (t+z, t+w)
=\langle \phi_\lambda (z), \phi_\lambda (w)\rangle .\eqno(3.12)$$
\indent We have, from the differential equation,
$$\Bigl( {{\partial }\over{\partial x}}+{{\partial }\over{\partial
y}}\Bigr) K_\lambda (x,y)=-\lambda \Bigl\langle {{E(x)-E(y)}\over{x-y}}
\Psi_\lambda (x), \Psi_\lambda (y)\Bigr\rangle  
-\Bigl\langle {{F(x)-F(y)}\over{x-y}}
\Psi_\lambda (x), \Psi_\lambda (y)\Bigr\rangle .\eqno(3.13)$$
\noindent By the hypotheses on $E$ and $F$, there exist
self-adjoint matrices $E_1, F_1\geq
0$, self-adjoint matrices $E_0$ and $F_0$, and positive matrix measures
$\Omega_E$ and $\Omega_F$ such that
$$E(x)=E_1x+E_0+\int_{\varepsilon}^\infty\Bigl( {{u}\over
{1+u^2}}-{{1}\over {u+x}}\Bigr)\Omega_E(du),\eqno(3.14)$$
\noindent and
$$F(x)=F_1x+F_0+\int_{\varepsilon}^\infty\Bigl( {{u}\over
{1+u^2}}-{{1}\over {u+x}}\Bigr)\Omega_F(du),\eqno(3.15)$$
\noindent hence (3.13) equals
$$-\lambda {{E(x)-E(y)}\over{x-y}}-{{F(x)-F(y)}\over{x-y}}$$
$$= -\lambda E_1-\lambda \int_{\varepsilon}^\infty {{\Omega_E (du)}\over
{(u+x)(u+y)}}-F_1-\int_{\varepsilon}^\infty {{\Omega_F (du)}\over
{(u+x)(u+y)}}.\eqno(3.16)$$
\noindent By a straightforward Hilbert space construction similar to
that in [3], we can
introduce $H_0$ and $\phi \in L^2((0, \infty );H_0)$ such that
$$\Bigl( {{\partial }\over{\partial x}}+{{\partial }\over{\partial
y}}\Bigr)
\Bigl( K_\lambda (x,y)\Bigr)=-\langle \phi_\lambda (x), \phi_\lambda
(y)\rangle_{H_0} .\eqno(3.17)$$
\noindent Hence we have
$$\Bigl( {{\partial }\over{\partial x}}+{{\partial }\over{\partial
y}}\Bigr)
\Bigl( K_\lambda (x,y)-\Gamma_{\phi_\lambda}^\dagger \Gamma_{\phi_\lambda}
(x,y)\Bigr)=0,\eqno(3.18)$$
\noindent and 
$$K_\lambda (x,y)-\Gamma_{\phi_\lambda}^\dagger \Gamma_{\phi_\lambda} (x,y)\rightarrow
0\qquad (x, y\rightarrow\infty)\eqno(3.19)$$
\noindent and hence 
$$K_\lambda (x,y)=\Gamma_{\phi_\lambda}^\dagger \Gamma_{\phi_\lambda}
(x,y)=\int_0^\infty \langle \phi_\lambda (s+x), 
\phi_\lambda (s+y)\rangle ds.\eqno(3.20)$$
\noindent Now $\Gamma_{\phi_\lambda}$ is bounded since $K_\lambda$ is bounded. Then for any Hankel operator  $S^\dagger_t\Gamma_{\phi_\lambda} =\Gamma_{\phi_\lambda} S_t$. Hence, $S^\dagger_t\Gamma_{\phi_\lambda}^\dagger \Gamma_{\phi_\lambda} S_t=\Gamma_{\phi_\lambda}^\dagger P_{(t, \infty )}\Gamma_{\phi_\lambda} $ so that
$$K_{t, \lambda }=S_t^\dagger K_\lambda S_t=
\Gamma_{\phi_\lambda}^\dagger P_{(t, \infty )} \Gamma_{\phi_\lambda}
\qquad (\lambda >0).\eqno(3.21)$$  
\indent (ii) When $H_0$ has finite dimension, the kernel $\langle
\phi_\lambda (x), \phi_\lambda (y)\rangle$ has finite rank by Lemma
3.2(ii).\par
\indent (iii) When $H_0={\bf C}$, the kernel of the Hankel operator is scalar-valued.\par
\indent (iv) In particular, the Hankel operator with $\phi_\lambda :(0, \infty )\rightarrow {\bf R}$ is self-adjoint and $K_\lambda =\Gamma^2_{\phi_\lambda}.$\par 
\rightline{$\square$}\par
\vskip.05in

\noindent {\bf Corollary 3.5.} {\sl Suppose that 
$K_\lambda=\Gamma_{\phi_\lambda}^{\dagger}\Gamma_{\phi_\lambda}$, as in
Theorem 3.4. Then under the Laplace transform ${\cal L}:L^2((0, \infty
); {\bf C}^{2m})\rightarrow H^2({\bf C}_+; {\bf C}^{2m})$, the nullspace of $K$ is
unitarily equivalent to $\theta H^2({\bf C}_+; {\bf C}^{2m})$ for some
bounded 
analytic function $\theta :{\bf C}_+\rightarrow M_{2m}({\bf C})$ 
that has unitary boundary values almost everywhere.}\par
\vskip.05in
\noindent {\bf Proof.} The null space of $K_\lambda$ equals the null space of
$\Gamma_{\phi_\lambda}$ and hence is a closed linear subspace of $L^2((0, \infty
); {\bf C}^{2m})$ which is invariant under the shift. Beurling's
theorem characterizes the images of such subspaces under the Laplace
transform; see [10].\par
\rightline{$\square$}\par
\vskip.05in
\indent {\sl Asymptotic forms of the differential equation as
$x,\lambda\rightarrow\infty$}\par
\indent We now consider circumstances under which (3.2) does
have a bounded or $L^2$ solution $\Psi_\lambda$. Suppose that 
$E$ and $F$ are is as in Theorem 3.4 and that 
$F_1=0$ in (3.14),
so that $F(x)$ is
bounded. Then there are the following basic cases (i), (ii) and (iii) for the
asymptotic form of (3.1) as $\lambda\rightarrow\infty$ and
$x\rightarrow\infty$.\par
\indent (i) Suppose that $E_1=0$. Then as $x\rightarrow\infty$ we have
$E(x)\rightarrow \tilde E_0$ where
$$\tilde E_0=E_0+\int_\varepsilon^\infty
{{u}\over{1+u^2}}\Omega_E(du),\eqno(3.22)$$
\noindent and the asymptotic form of the differential equation is 
$$J{{d}\over{dx}}\Phi_\lambda (x)=\lambda \tilde E_0\Phi_\lambda
(x)\eqno(3.23)$$
\noindent with solution
$$\Phi_\lambda (x) =\exp (-\lambda xJ\tilde E_0) \Phi_\lambda
(0).\eqno(3.24)$$
\noindent Now $\Re (J\tilde E_0)=[J, \tilde E_0]/2$, so 
$\Re (J\tilde E_0)$ is self-adjoint and has trace zero; hence 
$\Re (J\tilde E_0)$ is either zero, or has both positive eigenspaces and
negative eigenspaces. So in the following sub-cases, the solution of (3.23) has either:\par
\indent  (i)(a) $\Phi_\lambda $ constant;\par
\indent (i)(b)
$\Phi_\lambda $ oscillating boundedly as $x\rightarrow\infty$; or\par
\indent (i)(c) exponentially decaying solutions and
exponentially growing solutions.\par
\indent In sub-cases (i)(b) and (i)(c) there exist bounded solutions
$\Phi_\lambda $ to (3.23) such that 
$$K_\lambda (x,y)={{\langle J\Phi_\lambda (x), \Phi_\lambda
(y)\rangle}\over {x-y}}\eqno(3.25)$$
\noindent gives a bounded linear operator on $L^2(0, \infty ).$\par
\indent (ii) Suppose that $E_1$ is strictly positive definite. Then the asymptotic form
of the differential equation is 
$$J{{d}\over{dx}}\Phi_\lambda (x)=\lambda xE_1\Phi_\lambda
(x)\eqno(3.26)$$
\noindent with solution
$$\Phi_\lambda (x) =\exp (-\lambda x^2E_1/2) \Phi_\lambda (0),\eqno(3.27)$$ 
\noindent and we have sub-cases (a), (b) and (c) analogous to those in (i) above.\par
\indent (iii) $E_1$ of rank $1, \dots ,2m-1$. This case includes
variants of Airy's equation.\par
\vskip.05in
 
\noindent {\bf Examples 3.6.} 
\noindent (i) Sonine considered the one-parameter families of functions
$Y_\nu$ that satisfy the system
$$\eqalignno{Y_{\nu -1}+Y_{\nu +1}&={{2\nu }\over{z}}Y_\nu ,\cr
Y_{\nu -1}-Y_{\nu +1}&=2Y_\nu ',&(3.28)\cr}$$
\noindent as in [21, p 82]; the Bessel functions $Y_\nu
={\hbox{J}}_\nu$ give solutions. One can transform the differential equation
for ${\hbox{J}}_\nu$ into the system
$${{d}\over {dx}}
\left[ \matrix{u\cr v\cr}\right]=
J\left[ \matrix{-1/x-(1-\nu^2)/4x^2& 0\cr 0&-1\cr}\right]
\left[ \matrix{u\cr v\cr}\right],\eqno(3.29)$$
\noindent which has solution $u(x)=\sqrt{x} {\hbox{J}}_\nu (2\sqrt{x}).$ This
system is matrix monotone when $\nu =1$.\par
\indent (ii) Theorem 3.4 applies to the Airy kernel (1.3), which describes the {\sl soft edge} of certain matrix ensembles.
Likewise, the Bessel kernel describes the hard edge; see [3, 4] for
details.\par
\vskip.05in
\noindent {\bf 4. Determinantal random point fields}\par
\vskip.05in
\indent In section 3 we showed how some important kernels factorize as $K=\Gamma^\dagger \Gamma.$ Here we consider the properties of $\det (I-\lambda \Gamma^\dagger\Gamma )$; so first we introduce and solve the Gelfand--Levitan integral equation.\par
\vskip.05in
\noindent {\bf Lemma 4.1.} {\sl Suppose that $-A:D(A)
\subseteq H\rightarrow H$ is a generator of a bounded $C_0$ semigroup, $B:H_0\rightarrow D(A)$, and $C:D(A)\rightarrow H_0$ are linear operators, where $H_0$ has finite dimension $m$, and let $\phi (x)=Ce^{-xA}B$. Suppose further that $K_C(A), K_{B^\dagger} (A^\dagger ) \leq 1$.\par
\noindent (i) Then the $m\times m$ matrix kernel
$$T_\lambda (x,y)=-\lambda Ce^{-xA}\bigl( I+\lambda R_x\bigr)^{-1}e^{-yA}B\qquad (0<x<y, \vert\lambda\vert <1)\eqno(4.1)$$
\noindent gives the unique solution of the integral equation}
$$T_\lambda (x,y)+\phi (x+y)+\lambda \int_x^\infty T_\lambda (x,z)\phi (z+y)\, dz=0\qquad (0<x<y),\eqno(4.2)$$
\noindent {\sl and the kernel $T_\lambda(x,y)$ satisfies
$${{\partial^2T_\lambda }\over{\partial x^2}}-{{\partial^2 T_\lambda}\over{\partial y^2}}-q(x)T_\lambda (x,y)=0\eqno(4.3)$$
\noindent where $q(x)=-2{{d}\over{dx} }T_\lambda(x,x)$.\par
\noindent (ii) Suppose further that $m=1$, and that $\Theta_x$ and $\Xi_x$ are Hilbert--Schmidt. Then the determinant satisfies}
$$T_\lambda (x,x)=
{{d}\over{dx}}\log\det (I+\lambda
\Gamma_{\phi_{(x)}})\qquad (x>0).\eqno(4.4)$$
\vskip.05in
\noindent {\bf Proof.} (i) First, we have 
$\Vert R_x\Vert =\Vert\Xi_x\Theta_x^\dagger\Vert\leq 1,$
\noindent so $I+\lambda R_x$ is invertible and $T_\lambda (x,y)$ is well defined. One 
checks the identity by substituting the given expression for 
$T_\lambda $ into the integral equation. Further, $\Vert\Gamma_{\phi_{(x)}}\Vert=\Vert\Theta_x^\dagger \Xi_x\Vert$, 
so $I-\lambda \Gamma_{\phi_{(x)}^\dagger}$ is invertible
so hence solutions to the Gelfand--Levitan 
integral equation are unique.\par
\indent One can differentiate the integral equation and integrate by parts to obtain
$${{\partial^2T_\lambda }\over{\partial x^2}}-{{\partial^2 T_\lambda }\over{\partial y^2}}+\lambda q(x)\phi (x+y)+\lambda \int_x^\infty 
\Bigl({{\partial^2T_\lambda}\over{\partial x^2}}-{{\partial^2 T_\lambda }\over{\partial z^2}}\Bigr)\phi (z+y)\, dz=0,\eqno(4.5)$$
\noindent so by uniqueness 
$${{\partial^2T_\lambda }\over{\partial x^2}}-{{\partial^2 T_\lambda}\over{\partial y^2}}=q(x)T_\lambda (x,y).\eqno(4.6)$$
\indent (ii) When $H_0={\bf C}$ the kernel takes values in ${\bf C}$. Here $R_x=\Xi_x \Theta_x^\dagger$ is trace class, and we can rearrange the traces and compute 
$$\eqalignno{T_\lambda (x,x)&
=-\lambda {\hbox{trace}}\bigl(Ce^{-xA}(I+\lambda R_x)^{-1}
e^{-xA}B\bigr)\cr
&=-\lambda {\hbox{trace}}\bigl(
(I+\lambda R_x)^{-1}e^{-xA}BCe^{-xA}\bigr)\cr
&={{d}\over{dx}}{\hbox{trace}}\log (I+\lambda R_x)\cr
&={{d}\over{dx}}\log\det (I+\lambda \Gamma_{\phi_{(x)}}),&(4.7)}$$
\noindent where the last step follows from Proposition 2.2.\par 
\rightline{$\square$}\par
\noindent Our first application is to the context of Theorem
1.2(i), where we consider determinantal random point fields
associated with the observability Gramian.\par
\vskip.05in
\noindent {\bf Theorem 4.2.} {\sl Suppose that $(4^o)$ $\Theta :L^2((0, \infty ); 
{\bf C})\rightarrow H$ defines a Hilbert--Schmidt operator, and that
$(2^o)$ the operator norm is $\Vert \Theta\Vert <1.$\par
\noindent (i) Then there exists a determinantal random point field on
$(0, \infty )$ such that $\nu (x,\infty )$ is the number of points in
$(x, \infty )$ and such that the generating function satisfies}
$${\bf E} z^{\nu (x, \infty )}=\det (I+(z-1)Q_x)\qquad (z\in {\bf C},
x>0).
\eqno(4.8)$$
\noindent {\sl (ii) Let $F$ be the cumulative distribution function
$$F(x)=\cases{ {\bf P}[\nu (x, \infty )=0]\qquad (x\geq 0)\cr
    0\qquad (x<0).\cr}\eqno(4.9)$$
\noindent Then}
$${F'(x)}/{F(x)}={\hbox{trace}}
((A+A^\dagger)Q_x(I-Q_x)^{-1})\qquad (x>0).\eqno(4.10)$$
\noindent {\sl (iii) In particular, if $A=A^\dagger$, then $\det (I+(z-1)\Gamma_{\phi_{(x)}})$ gives a generating function.}\par
\noindent {\sl (iv) When $A=A^\dagger$, the kernels  
$$T_\lambda (x,y)=-\lambda Ce^{-xA}(I+\lambda Q_x)^{-1}
e^{-yA}C^\dagger\qquad (0<x<y, \vert \lambda\vert <1)\eqno(4.11)$$
\noindent and $\phi (x+y)=Ce^{-(x+y)A}C^\dagger$ satisfy the Gelfand--Levitan integral equation
$$T_\lambda (x,y)+\lambda \phi (x+y)+
\lambda\int_x^\infty T_\lambda (x,z)\phi (z+y)\, dz=0\eqno(4.12)$$
\noindent and the diagonal satisfies} 
$$T_\lambda (x,x)={{d}\over{dx}}\log\det (I+\lambda 
\Gamma_{\phi_{(x)}}).\eqno(4.13)$$
\vskip.05in

\vskip.05in
\noindent {\bf Proof.} (i) The kernel $Ce^{-sA}e^{-tA^\dagger}C^\dagger$
of $\Theta^\dagger\Theta$ gives an integral operator on $L^2(0,
\infty)$ such that $0\leq \Theta^\dagger \Theta\leq I$; hence by 
Lemma 1.1 is associated with a determinantal random point field such that 
$$\eqalignno{{\bf E} z^{\nu (x, \infty )}&=
\det (I+(z-1)\Theta^\dagger\Theta P_{(x, \infty )})\cr
&=\det (I+(z-1)\Theta P_{(x, \infty )}\Theta^\dagger )\cr
&= \det (I+(z-1)Q_x),&(4.14)}$$
\noindent so the determinant involving the observability Gramian gives rise to the determinantal random point field.\par 
\indent (ii) We consider the probability that all of the random points lie in $(0, x)$. The operator $I-Q_x$ is invertible since $\Vert Q_x\Vert_{op}<1$, and we have 
$$F(x)=\det (I-Q_x)\qquad (x>0).\eqno(4.15)$$
\noindent By a familiar formula for determinants, we have
$$\eqalignno{{{d}\over{dx}}\log \det
(I-Q_x)&={{d}\over{dx}}{\hbox{trace}}\log (I-Q_x)\cr
&=-{\hbox{trace}}\Bigl( (I-Q_x)^{-1}{{d}\over{dx}}Q_x\Bigr)\cr
&={\hbox{trace}}\Bigl( (I-Q_x)^{-1}(A^\dagger Q_x+Q_xA)\Bigr)\geq 0,
&(4.16)}$$
\noindent where the last step follows from the Lyapunov equation
$(1^o)$. 
Condition $(3^o)$ reassures us that $F(x)$ is indeed an increasing function, 
and that the associated probability density function satisfies
(4.10).\par
\indent (iii) If $A=A^\dagger$, then $Ce^{-(s+t)A}C^\dagger$ is the kernel of $\Gamma_{\phi_{(x)}}=\Theta_x^\dagger \Theta_x$, so 
$$\det (I+(z-1)Q_x)=\det (I+(z-1)\Gamma_{\phi_{(x)}}).\eqno(4.17)$$
\indent (iv) This is a special case of Lemma 4.1.\par
\rightline{$\square$}\par
\vskip.05in
\noindent Now we state the variant which arises in random matrix theory
as in Theorem 3.4(iv) and Theorem 1.2(ii).\par
\vskip.05in

\noindent {\bf Theorem 4.3.} {\sl Suppose that $A,B$ and $C$ 
satisfy the hypotheses of Lemma 4.1, and that $\phi =\phi^\dagger$. Let $\phi_{(x)}(y)=\phi (2x+y)$.\par
\noindent (i) Then there exists a determinantal random point field on
$(0, \infty )$ such that $\nu (x,\infty )$ is the number of points in
$(x, \infty )$ such that the generating function satisfies
$${\bf E} z^{\nu (x, \infty )}=
\det (I+(z-1)\Gamma_{\phi_{(x)}}^2)\qquad (z\in {\bf C},
x>0).\eqno(4.18)$$
\noindent (ii) Further, 
$${{d}\over{dx}}\log\det (I-\lambda^2 \Gamma_{\phi_{(x)}}^2)=
T_{-\lambda }(x,x)+T_{\lambda }(x,x)\qquad (\vert\lambda\vert
<1),\eqno(4.19)$$
\noindent where $T_\lambda (x,y)=-\lambda Ce^{-xA}(I+\lambda R_x)^{-1}
e^{-yA}B$ satisfies a Gelfand--Levitan equation as in (4.2).}\par
\vskip.05in
\noindent {\bf Proof.} First we check that $K_x=\Gamma_{\phi_{(x)}}^2$ satisfies $(2^o),(4^0)$ and $(5^0)$. We have $\Vert\Gamma_{\phi_{(x)}}\Vert\leq 1,$ so $0\leq K_x\leq I$, and $K_{x}=\Theta_x^\dagger R_x$ is of trace class. Hence by Soshnikov's theorem, we can form a determinantal random point field with generating function as above.\par
\indent To calculate the determinant, one can use the identity
$$\eqalignno{\log \det (I-\lambda^2K_x)&=\log\det 
(I-\lambda \Gamma_{\phi_{(x)}}) 
+\log\det (I+\lambda \Gamma_{\phi_{(x)}}),\cr
&=\log\det (I-\lambda R_x) +\log\det (I+\lambda R_x).&(4.20)}$$
\indent (ii) The terms on the right-hand side satisfy 
$${{d}\over{dx}}\bigl( \log\det (I-\lambda R_x)
+\log\det (I+\lambda R_x)\bigr)=T_{-\lambda} (x,x)+T_\lambda (x,x).$$
\noindent by the Gelfand--Levitan equation. Indeed we can replace $B$ 
in Lemma 4.1 by $\pm B$ to introduce $\pm \lambda\Gamma_{\phi_{(x)}}$.\par 
\rightline{$\square$}\par
\vskip.05in
\indent We defer discussion of Theorem 1.2(iii) until section 6. In
section 5, we consider the determinant in Theorem 4.3 from the
perspective of scattering theory.\par
\vskip.05in
\noindent {\bf 5. Scattering and inverse scattering}\par
\vskip.05in
\indent The Gelfand--Levitan integral equation of Lemma 4.1 
is closely connected to the Schr\"odinger equation, as we discuss in
this section. Our aim is to identify a group of bounded linear
operators which acts naturally on the $\phi$ that appear in Theorem
4.3, and hence on the determinants.\par
\indent Given $\phi (x)=Ce^{-xA}B$ as in Lemma 4.1, we can solve the Gelfand--Levitan 
equation and recover $q(x)=-2{{d}\over{dx}}T_\lambda (x,x).$ 
Further, given $T_\lambda$ as in Lemma 4.1, the function  
$$\psi (x;k)=e^{ikx}+\int_x^\infty e^{iyk}T(x,y)dy\eqno(5.1)$$
\noindent satisfies
$$-{{d^2}\over{dx^2}}\psi (x;k)+q(x)\psi (x;k)=k^2\psi (x;k)\eqno(5.2)$$
\noindent and
$$\psi (x;k)\asymp e^{ikx}\qquad (x\rightarrow\infty ).\eqno(5.3)$$
\noindent This is a straightforward calculation, based upon (2.25).\par

\indent (i) Let $(\lambda_j)_{j=1}^n$ be the discrete spectrum of $-{{d^2}\over{dx^2}}+q$ in $L^2({\bf R})$, written $\lambda_j=-\kappa_j^2$ with $\kappa_j>0$
so that each $\lambda_j=-\kappa_j^2$ is associated with an
eigenfunction $\psi (x; \lambda_j)$ that is asymptotic to $e^{-\kappa_jx}$ as 
$x\rightarrow\infty$ and $\kappa_n\geq \dots \geq \kappa_1>0$. We
take $c(-\kappa_j^2)$ to be a constant associated with
$-\kappa_j^2$.\par 
\indent (ii) The continuous spectrum is $\Sigma_c=[0, \infty )$, which has multiplicity two. For $k\in {\bf R}$ and $\lambda =k^2>0$, there 
exists a solutions $\psi (x;k)$ to (5.2) with asymptotic behaviour 
$$\psi (x,k)\asymp \cases{ e^{-ikx}+b(k)e^{ikx}&
as $x\rightarrow\infty $;\cr a(k)e^{-ikx}& as
$x\rightarrow -\infty$.\cr}\eqno(5.4)$$
\indent By [12], the reflection coefficient $b$ 
belongs to $C_0^\infty$, satisfies $b(0)=-1$ and $b(-k)=\bar b(k)$ and 
$$\int_{-\infty}^\infty \Bigl\vert{{\log (1-\vert b(k)
\vert^2)}\over{1+k^2}}\Bigr\vert dk <\infty .\eqno(5.5)$$
\par
\indent (iii) By results from [12, 13], the transmission coefficient $a$ extends to define the outer function 
$$a(k)=\exp\Bigl({{1}\over{2\pi i}}\int_{-\infty}^\infty
 {{\log  (1-\vert b(t)\vert^2)}\over {t-k}}\, dt\Bigr)\eqno(5.6)$$
on $\{ k: \Im k>0\}$ such that $\vert a(k)\vert=(1-\vert
b(k)\vert^2)^{1/2}$ for $k\in {\bf R}$.\par 
\vskip.1in
\noindent The scattering map $q\mapsto \phi $ associates, to the potential $q$, the function
$$\phi (x) =\sum_{j=1}^n c(-\kappa_j^2)^2e^{-\kappa_jx}+
{{1}\over{2\pi}}\int_{-\infty}^\infty b(k )e^{ikx} \, dk .\eqno(5.7)$$
\noindent where the eigenvalues $\lambda_j=-\kappa_j^2$, the
normalizing constants $c(-\kappa^2_j)$ and the reflection
coefficient $b(\kappa )$ are the scattering data. By (ii), $\phi (x)$
is real.\par
\indent The aim of inverse scattering is to recover $q$, up to translation from the scattering data. 
\noindent The following section uses computations which are extracted
 from [1,2], and originate in other calculations of inverse problems
as in [6,21]. 
The reflection coefficient $b$, the negative eigenvalues $\lambda_j$ and the normalising constants $c(-\kappa^2_j)$ determine $q$ uniquely up to translation.\par 
 
\indent  Our general approach to the inverse spectral problem is to go \par
$$\phi\mapsto (-A,B,C)\mapsto T\mapsto q.\eqno(5.8)$$
We now consider the first step in the process, namely realising $(-A,B,C)$ from a given $\phi$.\par
\vskip.05in 
\noindent {\bf Definition} {\sl (Realisation).} Given an bounded 
linear operator $\Gamma$, we wish to find a linear system
$(-A,B,C)$ such that the
 corresponding Hankel operator $\Gamma_\phi$ is unitarily equivalent 
to $\Gamma$. In particular, given scattering data $\phi (x)$, 
we wish to find a balanced linear system such that $\phi (x)=Ce^{-xA}B$ for $x>0$.\par
\vskip.05in
\indent The following Lemma gives a characterization up to 
unitary equivalence of a special class of Hankel operators.\par
\vskip.05in
\noindent {\bf Definition} {\sl (Spectral multiplicity).} 
For a self-adjoint and bounded linear operator $A$ on $H$ with spectrum
$S$, let 
$$H=\int_S^\oplus H(\lambda )\,\mu (d\lambda )$$
\noindent be the spectral resolution, where $\mu$ is a bounded
positive Radon measure on $S$, such that $Af(\lambda )=\lambda f(\lambda )$. 
 Now let $\delta (\lambda )={\hbox{dim}}H(\lambda )$ be the 
spectral multiplicity function for $\lambda\in S$.\par
\vskip.05in

\noindent {\bf Lemma 5.1.} {\sl Let $\Gamma$ be a self-adjoint and bounded linear operator on $H$ such that:\par
\noindent (i) the nullspace of $\Gamma$ is zero or infinite-dimensional;\par
\noindent (ii) $\Gamma$ is not invertible;\par
\noindent (ii) $\vert \delta (\lambda )-\delta (-\lambda )\vert\leq 1$ 
for $\mu$ almost all $\lambda $.\par
\noindent Then there exists a balanced linear system $(-A, B,C)$ with 
$H_0={\bf C}$ such that the Hankel operator $\Gamma_\phi$ on $L^2(0, \infty )$ with kernel $\phi (x+y)=Ce^{-(x+y)A}B$ is unitarily equivalent to $\Gamma$.}\par
\vskip.05in
\noindent {\bf Proof.} This is a special case of 
Theorem 1.1 on p. 257 of [14].\par
\rightline{$\square$}\par
\vskip.05in
\noindent {\bf Proposition 5.2.} {\sl Suppose that
 the hypotheses of Lemma 5.1 hold.\par
\noindent (i) If the spectral density function $b$ on the continuous spectrum is 
identically zero, then the scattering data can be 
realized by a linear system with finite dimensional $H$ and $H_0$.}\par
\noindent {\sl (ii) Suppose that in the corresponding linear system
 $A$ is a finite matrix such that all its eigenvalues $\kappa_j$ 
satisfy $\Re \kappa_j>0$. Then the system is admissible.}\par  
\vskip.05in
\noindent {\bf Proof.} (i) By a theorem of Fuhrmann [15], one can choose 
$A$ to be a finite-rank operator if and only if the transfer function $\hat \phi$ is a rational
 function which is analytic on the closure of ${\bf C}_+
\cup \{\infty \}$.\par
\indent (ii) The system is admissible since $\Vert
Ce^{-tA}\xi\Vert\leq Me^{-\kappa t}\Vert \xi\Vert$ for some $M,\kappa
>0$ and $\xi\in H_0$.\par
\rightline{$\square$}\par
\vskip.05in
\noindent {\bf Remarks.} (i) Not all self-adjoint Hankel operators 
satisfy the condition (ii) of Lemma 5.1. Consequently, there is a 
distinction between those self-adjoint 
Hankel operators that can be realised by linear systems in continuous 
time with one-dimensional input and output spaces and the more general
 class that can be realised by linear systems in discrete time.
 In this paper we concentrate on the continuous time case, while
 McCafferty has considered analogous results in discrete time, as
in [11].\par
\indent (ii) If $b=0$ in (5.9), then Propositions 5.2 and 2.3 
apply to $\phi$.\par
\vskip.05in
\noindent {\bf Definition} {\sl (Evolution)}. For a system $(-A,B,C)$, we refer to $\phi$ as scattering
data. Given a $C_0$ group $E(t)$ on $H$ and $D(A)$, we can form the system
$(-A, B, CE(t))$ and introduce $\tilde E(t)\phi
(s)=CE(t)e^{-sA}B$; thus
the scattering data evolves with $t$.\par
\vskip.05in
\indent Article [3] highlighted the importance of groups that
satisfy the Weyl relations;
here we show that these are associated with evolutions that do not
change the determinants in Theorem 1.2.\par
\vskip.05in 

\noindent {\bf Proposition 5.3.} {\sl Suppose that $D$ is skew
self-adjoint and generates a $C_0$ group of unitary operators
that satisfies
$e^{-sD}e^{-tA}=e^{-i\alpha st}e^{-tA}e^{-sD}$ for all $s\in {\bf R},
t\geq 0$ and some $\alpha\in {\bf R}$. Then $(-A, e^{sD}B,
Ce^{-sD})$ has observability Gramian $Q^{(s)}$, operator
$R^{(s)}$ and controllability 
Gramian $L^{(s)}$ such that (i) $\det
(I+\lambda Q^{(s)})$, (ii) $\det (I+\lambda R^{(s)})$ and (iii) 
$\det (I+\lambda Q^{(s)}L^{(s)})$ are independent of $s$.}\par 

\vskip.05in
\noindent {\bf Proof.} (i) We have $e^{-tA^\dagger}e^{sD}C^\dagger
Ce^{-sD}e^{-tA}=e^{sD}e^{-tA^\dagger}C^\dagger Ce^{-tA}e^{-sD}$, so\par
\noindent 
$Q^{(s)}=e^{sD}Q_0e^{-sD}$ and by unitary equivalence $\det (I+\lambda
Q^{(s)})=\det (I+\lambda Q_0).$\par
\indent (ii) We have $R^{(s)}=e^{-sD}R_{0}e^{sD}$, so $\det (I+\lambda
R^{(s)})=\det (I+\lambda R_0)$.\par 
\indent (iii) Likewise, $L^{(s)}=e^{sD}L_0e^{-sD}$, and there is a unitary
equivalence $Q^{(s)}L^{(s)}=e^{sD}Q_0L_0e^{-sD}$ leading to $\det
(I+\lambda Q^{(s)}L^{(s)})=\det (I+\lambda Q_0L_0)$.\par
\rightline{$\square$}\par
\vskip.05in
\noindent {\bf Lemma 5.4.} {\sl Suppose that $\Psi_t(x ;k)$ gives a
  differentiable family of vectors in ${\bf C}^{2m}$ such that
$${{d}\over{dt}}\Psi_t(x;k) =Z_t(x;k)\Psi_t(x;k),\eqno(5.9)$$
\noindent where 
$$Z_t(x;k)=\left[\matrix{
\alpha_t(x;k)&\beta_t(x;k)\cr -
\gamma_t(x;k)&-\alpha_t(x;k)\cr}\right]\qquad 
\eqno(5.10)$$
\noindent and $\alpha_t (x;k), \beta_t(x;k)$ and $\gamma_t (x;k)$ are 
 symmetric $m\times m$ matrices. \par
\noindent (i) Then with the bilinear form on ${\bf
C}^{2m}$, the family of kernels
$$K_{t,x}(\kappa ,k)={{\langle J\Psi_t(x;\kappa ), \Psi_t(x;k)
\rangle}\over{\kappa -k}}\qquad (k, \kappa\in {\bf R}, k\neq \kappa )$$
\noindent satisfies $(10^o)$ and}
$${{\partial}\over{\partial t}}K_{t,x}(\kappa ,k)=
\Bigl\langle J\Bigl( {{Z_t(x;\kappa )-Z_t(x;k)}\over{\kappa -k}}\Bigr)\Psi_t(x;\kappa ), 
\Psi_t(x;k)\Bigr\rangle\eqno(5.11)$$
\noindent {\sl (ii) If the $\alpha_t(x;k),
\beta_t(x;k)$
and $\gamma_t(x;k)$ are rational functions of $k$, then ${{\partial}\over{\partial t}} K_{t,x}$ is of
finite rank.}\par 
\vskip.05in
\noindent {\bf Proof.} (i) This follows by direct calculation, where 
the the effect of ${{\partial}\over{\partial t}}$ is to replace 
$J$ by $JZ_t(x;k)+Z_t(y;k)^TJ$ in the kernel. Then one uses the identities
$$JZ_t(x;\kappa )+Z_t(x;k)^T J=J(Z_t(x;\kappa )-Z_t(x;k)),\eqno(5.12)$$
\noindent which follow from the special form of the matrices. \par
\indent (ii) Given the formula (5.21), one can use the partial
fraction decomposition of the entries to express ${{\partial}\over{\partial t}}K_{t,x}(\kappa ,k)$ as a sum
of products of functions in the variable $\kappa $ or $k$.\par
\rightline{$\square$}\par
\vskip.05in
\indent In accordance with the approach of [7], we are particularly
interested in the case where $k\mapsto Z_t(x;k)$ is a polynomial
such that the leading coefficient has trace zero. For Schr\"odinger's 
equation, we can introduce such families of matrices associated with the KdV flow. 
Let $C_0^\infty$ be the space of functions $f:{\bf R}\rightarrow {\bf C}$ that are infinitely differentiable and such that $\vert x\vert ^j\vert f^{(\ell )}(x)\vert \rightarrow 0$ as $x\rightarrow \infty $ for $j,\ell =0, 1,\dots $.\par
\vskip.05in
\indent Suppose that $q$ satisfies that $q=v'+v^2$. 
Given a $\psi$ that satisfies $-\psi''+q\psi=k^2\psi$, 
we have a solution of the symmetric Hamiltonian system
$${{d}\over{dx}}\left[ \matrix{ \psi\cr \rho\cr}\right]= 
\left[ \matrix{ v&ik\cr ik&-v\cr}\right] 
\left[ \matrix{ \psi\cr \rho\cr}\right] .\eqno(5.13)$$ 
\noindent Now let $v$ evolve according to the modified Korteweg--de Vries equation
$$4{{\partial v}\over{\partial t}}={{\partial^3v}\over{\partial x^3}}
-6v^2{{\partial v}\over{\partial x}},\eqno(5.14)$$
\noindent and introduce functions of $(x,t)$ by
$$\alpha =(1/4)v_{xx}-(1/2)v^3, \quad \beta =(-1/2)(v_x+v^2),\quad
\gamma=
(1/2)(v_x-v^2),\quad \delta =v.\eqno(5.15)$$
\vskip.05in
\noindent {\bf Lemma 5.5.} {\sl The matrices
$$V_{t}(x;z)=\left[\matrix{ v&z\cr z&-v\cr}\right]\quad 
{\hbox{ and}}\quad Z_{t}(x;z)=
\left[\matrix{ \alpha +\delta z^2& \beta z+z^3\cr \gamma z+z^3&
-\alpha -\delta z^2\cr}\right]\eqno(5.16)$$
\noindent give a consistent system}
$$\cases{{{d}\over{dx}}\Psi =V_{t}(x;z)\Psi ,\cr
{{d}\over{dt}}\Psi =Z_{t}(x;z)\Psi .\cr}\eqno(5.17)$$
\vskip.05in
\noindent {\bf Proof.} As in [7], it follows by direct computation that 
$${{\partial V_{t}(x;z)}\over{\partial t}}-
{{\partial Z_{t}(x;z)}\over{\partial
x}}+[V_{t}(x;z),Z_{t}(x;z)]=0,\eqno(5.18)$$
\noindent so ${{\partial^2}\over{\partial x\partial
t}}\Psi={{\partial^2}\over{\partial t\partial x}}\Psi$ and 
the system is consistent. The key idea is that one can equate 
coefficients of the ascending powers of $z$, then one 
can eliminate the functions $\alpha ,\beta,\gamma$ and $\delta$ by simple calculus. \par
\rightline{$\square$}\par
\vskip.05in
\indent Let $\Psi_t(x;k)$ be the solution of (5.17) that
corresponds to $z=ik$ where $k$ belongs to ${\bf R}$ and $k^2$ to the continuous spectrum
$(0,\infty )$. With the bilinear form $\langle .\, ,\,.\,\rangle$ on ${\bf
C}^2$, let
$$K_{t,x}(\kappa ,k)={{\langle J\Psi_t(x;\kappa
),\Psi_t(x;k)\rangle}\over {i(\kappa -k)}}.\eqno(5.19)$$
\noindent where the numerator vanishes on $k=\kappa$, so $K_{t,x}$
is an integrable operator. This family of operators undergoes a
natural evolution under the KdV flow, as follows.\par
\vskip.05in
\noindent {\bf Theorem 5.6.} {\sl Suppose that $\Psi_t(x;k)$ give a 
locally bounded family of solutions which is differentiable in
$(t,x,k)$ and subject to $\Psi_t(0;k)=\Psi_t$ for some $\Psi_t\in
{\bf C}^2$.\par
\noindent (i) If $v(x)=0$, then $K_{t,x}$ is a multiple of the sine
kernel [3, (1.2)].\par
\noindent (ii) The kernels $K_{t,x}$ satisfy $(10^o)$ and $(5^o)$, so
${{\partial}\over{\partial x}}K_{t,x}$ and ${{\partial}\over{\partial
t}}K_{t,x}$ are of finite rank.\par
\noindent (iii) The function $u={{\partial v}\over{\partial x}}+v^2$
satisfies the KdV equation, and as $q(x)$ evolves to $u(x,t)$ the scattering data for $u$ undergoes a linear evolution
$\phi \mapsto E(t)\phi$.\par
\noindent (iv) Let $(b, c(-\kappa_j^2), \kappa_j)$ be the scattering data for $q(x)$, and suppose that 
$b(k),b'(k)$ and $k^2b(k)$ belong to
$L^2({\bf R}; dk)$. Then $\Gamma_{E(t)\phi}$ gives a
Hilbert--Schmidt operator for all $t$.}\par

\vskip.05in
\noindent {\bf Proof.} (i) This is an elementary computation.\par
\indent (ii) Using Lemma 5.4, we calculate the derivatives, and find
$${{\partial K_{t,x}}\over{\partial t}}=\Bigl\langle \left[\matrix {
-\gamma +(k^2+k\kappa +\kappa^2)& i\delta (k+\kappa )\cr 
      i\delta (k+\kappa )& \beta -(k^2+k\kappa +k^2)\cr}\right] 
\Psi_t(x;\kappa
), \Psi_t(x;k)\Bigr\rangle ,\eqno(5.20)$$
\noindent which gives a kernel of finite rank, and likewise
$${{\partial K_{t,x}}\over{\partial
x}}=\Bigl\langle\left[\matrix{1&0\cr 0&-1\cr}\right] \Psi_t(x;\kappa
), \Psi_t(x,k)\Bigr\rangle \eqno(5.21)$$
\noindent which also gives a kernel of finite rank on $L^2(-\infty, \infty )$.\par 
\indent (iii) By Miura's transformation, the function $u={{\partial v}\over {\partial x}}+v^2$ satisfies the KdV equation
$$4{{\partial u}\over{\partial t}}={{\partial^3 u}\over{\partial x^3}}
-6u{{\partial u}\over{\partial x}};\eqno(5.22)$$
\noindent see [6, p.65]. The evolution of the potentials
$u(x,0)\mapsto u(x,t)$ under the KdV flow gives rise to a linear
evolution on the scattering data.
 Now let $\Psi_{t}(x;k)$ be a continuous and uniformly bounded family of solutions of the system
$$\cases{{{d}\over{dx}}\Psi_{t}(x;k) =U_{t}(x;k)\Psi_t(x;k)\cr
{{d}\over{dt}}\Psi_t(x;k) =W_{t}(x;k)\Psi_{t}(x;k)\cr}\eqno(5.23)$$
\noindent where 
$$U_{t}(x;k)=\left[\matrix{ 0&1\cr u-k^2 &0\cr}\right],\quad
W_{t}(x;k)={{-1}\over{4}}\left[\matrix{ 
4ik^3-{{\partial u}\over{\partial x}}&2u+4k^2\cr 
2(u+2k^2)(u-k^2)-{{\partial^2u}\over{\partial x^2}}&4ik^3+{{\partial u}\over {\partial x}}\cr}\right].\qquad 
\eqno(5.24)$$
\noindent Then by considering the shape of the matrices in (5.24), we
obtain the asymptotic forms of the solutions
$$\Psi_t(x;k)\asymp -ia(k)e^{-ikx}
\left[\matrix{i\cr k\cr}\right]\qquad (x\rightarrow -\infty ),$$
$$\Psi_t(x;k)\asymp -ie^{-ikx}\left[\matrix{i\cr k\cr}\right]
+ib(k)\left[\matrix{-i\cr k\cr}\right]e^{ikx-2ik^3t}\qquad (x\rightarrow 
\infty );\eqno(5.25)$$
\noindent hence $a(k)\mapsto a(k)$ and $b(k)\mapsto b(k)e^{-2ik^3t}$
under the flow. By [6, p. 75], there is a group of linear operators $E(t)$ on the Hilbert space ${\bf C}^n\oplus L^2({\bf R})$ defined by
$$E(t)\phi (x)=\sum_{j=1}^n c(-\kappa_j^2)^2 
e^{-\kappa_jx-2\kappa_j^3t}
 +{{1}\over{2\pi}}\int_{-\infty}^\infty b(k) e^{ikx-2ik^3t}
dk\eqno(5.26)$$
\noindent such that $u(x,y)$ corresponds to $E(t)\phi$, and 
$\Vert E(t)\Vert=\max\{ e^{-2t\kappa^3_n}, 1\}.$\par
\noindent By applying Fourier inversion to the definition of the Airy
function, we can express the integral over the continuous spectrum
as
$$ {{1}\over{2\pi}}\int_{-\infty}^\infty b(k)
e^{ikx-2ik^3t}dk={{-1}\over{(6t)^{1/3}}}\int_0^\infty
{\hbox{Ai}}\Bigl( {{x-y}\over{(-6t)^{1/3}}}\Bigr)\phi (y)dy.$$
\indent (iv) The Hankel operator with kernel $\sum_{j=1}^n
c(-\kappa_j^2)^2 e^{-\kappa_j(x+y)-2\kappa_j^3t}$ is clearly of 
trace class, so we need to consider the Hankel arising from $b$. To show that $\int_0^\infty x(E(t)\phi (x))^2dx<\infty$, it
suffices by Plancherel's theorem to show that $b(k)$ and  
${{d}\over{dk}}(b(k)e^{-2ik^3t})$ belong to $L^2({\bf R}; dk).$ This
follows directly from the hypotheses.\par
\rightline{$\square$}\par
\vskip.05in
\noindent {\bf 6 Determinantal random point fields associated with the Zakharov--Shabat system}\par
\vskip.05in
\noindent In this final section we prove the remaining case (iii) of
Theorem 1.2, and then we address the corresponding scattering
theory.\par
\noindent Consider the matricial Gelfand--Levitan integral equation
$$T(x,y)+\lambda \Phi (x+y)+\lambda \int_x^\infty T(x,z)\Phi (z+y)\,
 dz=0\qquad (0<x<y),\eqno(6.1)$$
\noindent where, suppressing the dependence of $T$ upon $\lambda$, we
write
$$T(x,y)=\left[\matrix{ \bar U(x,y)& V(x,y)\cr -\bar 
V(x,y)&U(x,y)\cr}\right] ,\qquad \Phi (x)=
\left[ \matrix{ 0& \bar \phi (x)\cr -\phi (x)&0\cr}\right].\eqno(6.2)$$ 
\vskip.05in
\noindent {\bf Theorem 6.1.} {\sl Suppose that the system 
$(-A, B,C)$ has $H_0={\bf C}$ and with $\phi_{(x)}(y)=Ce^{-(2x+y)A}B$ 
satisfies, as in Lemma 4.1:\par
\noindent $(2^o)$ $\Vert \Theta_x\Vert <1$ and $\Vert\Xi_x\Vert\leq 1$, and\par
\noindent $(4^o)$ $\Theta_x$ and $\Xi_x$ are Hilbert--Schmidt. \par
\noindent i) Then there exists a determinantal random point field on
$(0, \infty )$ such that $\nu (x,\infty )$ is the number of points in
$(x, \infty )$ and such that the generating function satisfies}
$$g_x(z)={\bf E} z^{\nu (x, \infty )} =\det (I+(z-1)\Gamma_{\phi_{(x)}}\Gamma_{\phi_{(x)}}^\dagger).\eqno(6.3)$$
\noindent {\sl (ii) Further ${{\partial }\over{\partial x}}
\log g_x(z)=2U(x,x)$, where $U$ is given by the diagonal of the 
solution of the Gelfand--Levitan equation (6.1).}\par
\vskip.05in
\noindent {\bf Proof.} (i) We have $0\leq Q_x\leq I$ and $0\leq L_x\leq I$.
Let $K_x=P_{(x, \infty )}\Theta^\dagger L_x\Theta P_{(x, \infty )}$, 
which defines a trace-class kernel on $L^2(x, \infty )$ and satisfies 
$0\leq K_x\leq I.$ Then by Lemma 1.1, there exists a determinantal random point field on $(x, \infty )$ with generating function
$$\eqalignno{\det (I+(z-1)K_x)&=\det (I+(z-1)
\Theta P_{(x, \infty )}\Theta^\dagger L_x)\cr
&=\det (I+(z-1)Q_xL_x),&(6.4)}$$
\noindent and we continue to rearrange this, obtaining
$$\eqalignno{\det (I+(z-1)K_x)&=\det (I+(z-1)P_{(x,\infty )}
\Theta^\dagger \Xi P_{(x, \infty )}\Xi^\dagger\Theta )\cr
&=\det (I+(z-1)P_{(x, \infty )}\Gamma_\phi 
P_{(x, \infty )}\Gamma_\phi^\dagger P_{(x, \infty )})\cr
&=\det (I+(z-1)\Gamma_{\phi_{(x)}}\Gamma_{\phi_{(x)}}^\dagger
).&(6.5)\cr}$$  
\indent (ii) This identity is proved in the following two Lemmas.\par
\vskip.05in 
\noindent {\bf Lemma 6.2.} {\sl Suppose that $H_0={\bf C}$ and let 
$\phi (x) =Ce^{-xA}B$ and $G_x=I+\lambda^2Q_xL_x$. Then the Gelfand--Levitan integral equation (6.1) reduces to
$$V(x,y)+\lambda \bar\phi (x+y)+
\lambda^2\int_x^\infty \!\!\!\int_x^\infty V(x,s)\phi (s+z)\bar\phi 
(y+z )dsdz=0\qquad (0<x<y),\eqno(6.6)$$
\noindent which has solution}
$$\eqalignno{V(x,y) &=-\lambda B^\dagger e^{-A^\dagger x}G_x^{-1}
e^{-A^\dagger y}C^\dagger \qquad (0<x<y),&(6.7)\cr
\bar U(x,y)&=\lambda \int_x^\infty V(x,z)\phi (z+y)\, dz.&(6.8)\cr}$$
\vskip.05in
\noindent {\bf Proof.} Once we have $V$, we can introduce $U$ via (6.7), and the resulting matrix $T$ satisfies the Gelfand--Levitan integral equation.
 To verify the equation for $T$, we first check that $G_x$ is invertible when $\Re \lambda^2>-1$. The operators $Q_x$ and $L_x$ are Hilbert--Schmidt and positive, so the operator $Q_xL_x$ is trace class, and hence the determinant satisfies
$$\det G_x=\det (I+\lambda^2 Q_x^{1/2}L_xQ_x^{1/2})>0\eqno(6.9)$$
\noindent since $\Re (I+\lambda^2Q_x^{1/2}L_xQ_x^{1/2})\geq (1-\Re \lambda^2)I.$\par
\indent 
One can postulate a solution of the form $V(x,y)=X(x)^\dagger e^{-A^\dagger y}C^\dagger ,$ for some function $X:(0, \infty )\rightarrow H$ and by substituting this into the integral equation, one finds that $X$ should satisfy
$$X(x)^\dagger e^{-A^\dagger y}C^\dagger
 \lambda B^\dagger e^{-A^\dagger (x+y)}C^\dagger$$
$$ +\lambda^2\int_x^\infty \int_x^\infty X(x)^\dagger 
e^{-A^\dagger s}C^\dagger Ce^{-A(s+z)}BB^\dagger e^{-A^\dagger (z+y)}
 C^\dagger dsdz=0,\eqno(6.10)$$
\noindent so we want 
$$X(x)^\dagger (I+\lambda^2Q_xL_x)+\lambda B^\dagger 
e^{-A^\dagger x}=0,\eqno(6.11)$$
\noindent and we can make this choice since $G_x$ is invertible.\par
\rightline{$\square$}\par
\vskip.05in
\noindent {\bf Lemma 6.3} {\sl The diagonal of the solution satisfies}
$$U(x,x)={{d }\over{d x}}{{1}\over{2}}\log\det
(I+\lambda^2\Gamma_{\phi_{(x)}}\Gamma_{\phi_{(x)}}^\dagger).\eqno(6.12)$$
\vskip.05in
\noindent {\bf Proof.} From (6.8), we have
$$\eqalignno{\bar U(x,y)&=-\lambda^2\int_x^\infty B^\dagger 
e^{-A^\dagger x}G_x^{-1}e^{-A^\dagger z} C^\dagger Ce^{-A(z+y)} Bdz\cr
&=-\lambda^2B^\dagger e^{-A^\dagger x}G_x^{-1}Q_xe^{-Ay}B.&(6.13)\cr}$$
 \noindent Hence we can write
$$\eqalignno{U(x,x)&=-\lambda^2B^\dagger e^{-A^\dagger x}
G_x^{-1} Q_xe^{-Ax}B\cr
&=-\lambda^2{\hbox{trace}}\Bigl( G_x^{-1} Q_x{{d}\over{d x}}L_x\Bigr)
.&(6.14)\cr}$$
\noindent We temporarily assume that $\lambda$ is real to 
derive certain identities, and then use analytic continuation to 
obtain them in general. Using Proposition 2.6, and rearranging 
various traces, we can derive the expressions
$$\lambda^2{\hbox{trace}}\Bigl( G_x^{-1} Q_x{{d L_x}
\over{dx}}\Bigr)={\hbox{trace}}\Bigl((G_x^\dagger )^{-1}A-A+(G_x)^{-1}A^\dagger -A^\dagger\Bigr)\eqno(6.15)$$ 
 \noindent and likewise 
$$\lambda^2{\hbox{trace}}\Bigl( G_x^{-1}{{dQ_x}
\over{dx}}L_x\Bigr)={\hbox{trace}}\Bigl((G_x^\dagger )^{-1}A-A+(G_x)^{-1}A^\dagger -A^\dagger\Bigr),\eqno(6.16)$$
\noindent and since ${{d G_x}\over {dx}}= 
\lambda^2({{dQ_x}\over {dx}}L_x+Q_x{{dL_x}\over {dx}}),$ we deduce that 
$$\eqalignno{U(x,x)&={{1}\over{2}}{\hbox{trace}}\Bigl( G_x^{-1}
{{dG_x}\over{dx}}\Bigr)\cr
&={{1}\over{2}}{{d}\over{dx}}
\log\det G_x\cr
&={{1}\over{2}}{{d}\over{dx}}\log\det (I+\lambda^2\Gamma_{\phi_{(x)}}
\Gamma_{\phi_{(x)}}^\dagger ).&(6.17)}$$
\noindent This concludes the proof of the Lemma, hence of
 Theorem 6.1(ii) and Theorem 1.2(iii).\par
\rightline{$\square$}\par
\vskip.05in

\indent We let $q\in C_0^\infty ({\bf R}; {\bf C})$ and consider 
the Zakharov--Shabat system 
$${{d}\over{dx}}\Psi (x;k) =\left[\matrix{-ik&q(x)\cr -\bar q(x)& ik\cr}
\right]\Psi (x;k) \eqno(6.18)$$
\noindent with $\Psi (x;k)$ a complex $2\times 2$ matrix. We observe 
that this matrix is skew-symmetric with zero trace, so the norm of
any solution is invariant under the evolution, as is the Wronskian of
any pair of solutions; hence the fundamental solution matrix of this 
system belongs to $SU(2)$. We introduce the solutions $\Psi_{+}(x;k), 
\Psi_{-}(x;k)\in SU(2)$ such that 
$$\Psi_{+}(x;k)\asymp \left[\matrix{ e^{-ikx}&0\cr 0&e^{ikx}\cr}
\right]\qquad (x\rightarrow\infty ),\eqno(6.19)$$
$$\Psi_{-}(x;k)\asymp \left[\matrix{ 
e^{-ikx}&0\cr 0&e^{ikx}\cr}\right]\qquad (x\rightarrow -\infty
);\eqno(6.20)$$
\noindent then we introduce the scattering matrix $S(k)\in SU(2)$ 
such that $\Psi_{-}(x;k)=\Psi_{+}(x;k)S(k)$ and we write
$$S(k)=\left[\matrix{ \alpha (k)&\hat\beta (k)\cr \beta (k)&-\hat\alpha (k)\cr}\right].\eqno(6.21)$$
\noindent Now suppose that $\alpha$ and $\beta$ are analytic 
on the upper half-plane, and that $\alpha$ has zeros at $\kappa_j$. 
As in [6], we introduce the scattering data
$$\phi (x) =\sum_{j=1}^n {{\beta (\kappa_j)}\over{\alpha' (\kappa_j)}}e^{i\kappa_jx}+{{1}\over{2\pi}}\int_{-\infty}^\infty {{\beta (k )}\over{\alpha (k)}}e^{ikx} \, dk .\eqno(6.22)$$
\noindent The sum contributes a function that decays exponentially as $x\rightarrow\infty$.\par
\vskip.05in
\noindent {\bf Proposition 6.4.} {\sl Let $(-A, B,C)$ realise the scattering data $\phi$ of the ZS system, suppose that the Gramians $Q_x$ and $L_x$ are Hilbert--Schmidt. Then the potential satisfies}
$$\vert q(x)\vert^2={{1}\over{2}}{{d^2}\over{dx^2}}
\log \det
(I+\lambda^2\Gamma_{\phi_{(x)}}\Gamma_{\phi_{(x)}}^\dagger ).\eqno(6.23)$$
\vskip.05in
\noindent {\bf Lemma 6.5.} {\sl Let $T$ be as in Lemma 6.2 and (6.2), and let
$$\Psi (x;k)=\left[ \matrix{ ae^{ikx}\cr be^{-ikx}\cr}\right]
 +\int_x^\infty T(x,y) \left[ \matrix {ae^{iky}\cr
 be^{-iky}\cr}\right]\, dy.\eqno(6.24)$$
\noindent Then 
$$-{{d^2}\over{dx^2}}
\Psi (x;k)+W(x)\Psi (x;k)=k^2\Psi (x;k)\eqno(6.25)$$
\noindent where} $W(x)=-2{{d}\over{dx}}T(x,x).$\par
\vskip.05in
\noindent {\bf Proof.} One can follow the proof of Lemma 4.1 and 
deduce that 
$${{\partial^2}\over{\partial x^2}}T(x,y) 
-{{\partial^2 }\over{\partial y^2}}T(x,y)=W(x)T(x,y).\eqno(6.26)$$
\noindent Then one can verify the differential equation for 
$\Psi (x;k)$ by direct calculation.\par
\indent From the original differential equation (6.18) we have 
$$-{{d^2}\over{dx^2}} \Psi (x;k)+
\left[\matrix{ -\vert q\vert^2&q'\cr-\bar q'&-\vert q
\vert^2\cr}\right]\Psi (x;k)=k^2\Psi (x;k),\eqno(6.27)$$
\noindent so by equating the matrix potential with $-2{{d}\over{dx}}T(x,x)$, we obtain
$$\left[\matrix{ -\vert q\vert^2&q'\cr-\bar q'&-\vert q\vert^2\cr}\right]=-2{{d}\over{dx}}
\left[\matrix{ \bar U(x,x)& V(x,x)\cr -\bar V(x,x) &U(x,x)\cr}\right].
\eqno(6.28)$$
\rightline{$\square$}\par
\vskip.05in
\indent Finally, we consider how the potential $q(x)$ evolves to
$u(x,t)$ under the
nonlinear Schr\"odinger equation. Suppose that 
$$W_t(x;\zeta )=\left[\matrix{ -i\zeta &u\cr -\bar
u&i\zeta\cr}\right],\qquad Z_t(x;\zeta )=\left[\matrix{-i\vert u\vert^2+2i\zeta^2& -i{{\partial
u}\over{\partial x}}-2u\zeta\cr 
-i\overline{{{\partial u}\over{\partial x}}} +2\bar u\zeta &i
\vert u\vert^2-2i\zeta^2\cr}\right].\eqno(6.29)$$
\vskip.05in
\noindent {\bf Proposition 6.6.} {\sl Suppose that $u$ satisfies the
nonlinear Schr\"odinger equation
$$i{{\partial u}\over{\partial t}}={{\partial^2u}\over{\partial
x^2}}+2\vert u\vert^2 u.\eqno(6.30)$$
\noindent (i) Then the pair of differential equations
$$\cases{{{d}\over{dx}}\Psi =W_t(x;\zeta )\Psi\cr
{{d}\over{dt}}\Psi =Z_t(x;\zeta )\Psi\cr}\eqno(6.31)$$
\noindent gives a consistent system.\par
\noindent (ii) Let $\Psi_t(x,\zeta )$ be family of solutions to (6.31)
with initial value $\Psi_t(0, \zeta )=\Psi_t$. Then the kernels
$$K_{t,x} (\kappa ,k)={{\langle J\Psi_t(x;\kappa ),
\Psi_t(x;k)\rangle}\over{\kappa -k}}\eqno(6.32)$$
\noindent satisfy $(10^o)$ and $(5^o)$; so ${{\partial}\over{\partial x}}K_{t,x}$
and ${{\partial }\over{\partial t}}K_{t,x}$ have finite rank.\par
\noindent (iii) The evolution of the potentials under the NLSE gives
rise to a linear evolution on the scattering data.}\par
\vskip.05in
\noindent {\bf Proof.} This is similar to that of Theorem 5.6.\par
\rightline{$\square$}\par
\vskip.05in
\noindent {\bf Acknowledgement.} {\smalletters I am grateful to Andrew McCafferty,
Stephen Power and Leonid Pastur for helpful conversations.}\par
 \vskip.1in 

\noindent {\bf References}\par
\noindent [1] M.J. Ablowitz and H. Segur, Solitons and the inverse
scattering transform, Society for Industrial and Applied Mathematics,
Philadelphia, 1981.\par
\noindent [2] T. Aktosun, F. Demontis, C van der Mee, Exact 
solutions to the focusing nonlinear Schr\"odinger equation, Inverse problems 23 (2007), 2171--2195.\par
\noindent [3] G. Blower, {Operators associated with the soft and hard edges from unitary
ensembles}, {J. Math. Anal. Appl.} {337} (2008), 239--265.
\par
\noindent [4] G. Blower, {Integrable operators and the squares
of Hankel operators, {J. Math. Anal. Appl.} 340 (2008), 943--953.\par

\noindent [5] S. Clark and F. Gesztesy, Weyl--Titchmarsh $M$-function
asymptotics for matrix-valued Schr\"odinger operators, {Proc. London
Math. Soc.} (3) 82 (2001), 701--724.\par
\noindent [6] P.G. Drazin and R.S. Johnson, Solitons: an introduction, 
Cambridge University Press, Cambridge, 1989.\par
\noindent [7] J.W. Helton, Operator theory, analytic functions, matrices, and electrical engineering, 
CBMS Regional conference series, American Mathematical Society,
1986.\par
\noindent [8] E. Hille, Lectures on ordinary differential equations, Addison Wesley, 1969.\par
\noindent [9] B. Jacob, J.R. Partington, and S. Pott, 
Admissible and weakly admissible observation operators, for the 
right shift semigroup, Proc. Edin. Math. Soc. (2) 45 (2002), 353--362.\par 
\noindent [10] P. Koosis, Introduction to $H_p$ spaces, Cambridge
University Press, Cambridge, 1980.\par
\noindent [11] A. McCafferty, Operators and special functions in random matrix theory, PhD Thesis, Lancaster 2008.\par
\noindent [12] H.P. McKean, The geometry of KdV II: three examples, J. Statist. Physics 46 (1987), 1115--1143.\par

\noindent [13] H.P. McKean, Geometry of KDV(3): determinants and unimodular isospectral flows, Comm. Pure Appl. Math. 45 (1992), 389--415.\par
\noindent [14] A.V. Megretski\u\i, V.V. Peller and S. Treil, The inverse spectral problem for self-adjoint Hankel operators, Acta Math. 174 (1995), 241--309.\par 
\noindent [15] N.K. Nikolskii, Operators, functions and systems; 
an easy reading. Volume 2: model operators and systems, 
(American Mathematical Society, 2002.\par
\noindent [16] A.G. Soshnikov. Determinantal random point fields, 2000, 
{\sl arXiv.org:math/0002099}\par
\noindent [17] C.A. Tracy and H. Widom, {Level-spacing distributions and the Airy
kernel,} {Comm. Math. Phys.} {159} (1994), 151--174.\par
\noindent [18] C.A. Tracy and H. Widom, {Level spacing distributions and the Bessel
kernel,} {Comm. Math. Phys.} {161} (1994), 289--309.\par
\noindent [19] C.A. Tracy and H. Widom, {Fredholm determinants, differential
equations and matrix models}, {Comm. Math. Phys.} {163} (1994), 33--72.\par
\noindent [20] C.A. Tracy and H. Widom, Correlation functions, cluster
functions and spacing distribution for random matrices, J. Statist.
Phys. 92 (1998), 809--835.\par

\noindent [21] E.T. Whittaker and G.N. Watson, A Course of Modern
Analysis, Fourth edition, Cambridge University Press, Cambridge, 1965.\par
\noindent [22] V.E. Zakharov and P.B. Shabat, A scheme for integrating the nonlinear equations of mathematical physics by the method of the inverse scattering problem. I, Funct. Anal. Appl. 8 (1974), 226--235.\par 
\vfill
\eject
\end